\documentclass[11pt]{article}
\usepackage{graphicx}
\usepackage{amsmath,amssymb,amsthm,amsfonts}
\usepackage{color}
\usepackage{amssymb}
\usepackage{cite}
\newtheorem{thm}{Theorem}[section]

\newtheorem{lem}[thm]{Lemma}
\newtheorem{prop}[thm]{Proposition}

\theoremstyle{definition}

\theoremstyle{remark}
\usepackage{appendix}
\newtheorem{rem}{Remark}[section]

\numberwithin{equation}{section}

\DeclareMathSymbol{\C}{\mathalpha}{AMSb}{"43}

\newcommand{\eps}{\varepsilon}

\newcommand{\lam}{\lambda}
\newcommand{\alp}{\alpha}

\newcommand{\ep}{\epsilon_a}

\newcommand{\R}{{\mathbb{R}}}
\newcommand{\h}{{\mathcal{H}}}
\newcommand{\inte}{\int_{\mathbb{R}^2}}

\def\e{\varepsilon_a}

\def\R{{\mathbb R}}
\def\C{{\mathbb C}}

\newcommand{\bsub}{\begin{subequations}}
\newcommand{\esub}{\end{subequations}$\!$}

\begin{document}
\title{The Nonexistence of Vortices for Rotating Bose-Einstein Condensates in Non-Radially Symmetric Traps}
\author{Yujin Guo\thanks{School of Mathematics and Statistics, and Hubei Key Laboratory of
Mathematical Sciences, Central China Normal University, P.O. Box 71010, Wuhan 430079, P. R. China. Email: \texttt{yguo@ccnu.edu.cn}. Y. J. Guo is partially supported by NSFC under Grant 11931012.
}
}

\date{\today}

\smallbreak \maketitle

\begin{abstract}
As a continuation of \cite{GLY}, we consider ground states of rotating Bose-Einstein condensates with attractive interactions in non-radially harmonic traps $V(x)=x_1^2+\Lambda ^2x_2^2 $, where $0<\Lambda \not =1$ and $x=(x_1, x_2)\in \R^2$.  For any fixed rotational velocity $0\le \Omega <\Omega ^*:=2\min \{1, \Lambda\}$, it is known that ground states exist if and only if $ a<a^*$ for some critical constant $0<a^*<\infty$, where $a>0$ denotes the product for the number of particles times the absolute value of the scattering length. We analyze the asymptotic expansions of ground states as $a\nearrow a^*$, which display the visible effect of $\Omega$ on ground states. As a byproduct, we further prove that ground states do not have any vortex in the region $R(a):=\{x\in\R^2:\,|x|\le C (a^*-a)^{-\frac{1}{12}}\}$ as $a\nearrow a^*$ for some constant $C>0$, which is independent of $0<a<a^*$.
\end{abstract}	
	

\vskip 0.05truein


\noindent {\it Keywords:} Bose-Einstein condensate; rotational velocity; nonexistence of vortices; limit profiles

\vskip 0.2truein



\section{Introduction}
Bose-Einstein condensate (BEC) is a state of matter, in which atoms or particles are cooled to  the sufficiently low temperature that a large fraction of them ``condense" into a single quantum state. The BECs in magneto-optical traps present remarkable phenomena, once the traps are set in the rotational motion. Actually, starting from the first physical achievement of rotating BECs in the late 1990s, various interesting quantum phenomena have been observed in the experiments of rotating BECs, including the critical-mass collapse \cite{Hulet1,D,HM}, the center-of-mass rotation \cite{Abo,LC,F}, and the appearance of quantized vortices \cite{A,CC,F}. Therefore, numerical simulations and mathematical theories of rotating BECs have been a focus of international interest in physics and mathematics over the past two decades, see \cite{Abo,A,Anderson,CR,D,F,IM-1,IM-2,Lieb06,LSS,LSY}. We also refer \cite{A,BB,CP,SSbook,S} for the mathematics of the relevant superfluids and superconductors.

The interactions between cold atoms in the condensates can be either repulsive or attractive, cf. \cite{A,B,D,F}.
For the repulsive case, the complex structures, including the quantized vortices, of the rotating trapped BECs were analyzed and simulated extensively in the past few years, see \cite{Abo,A,AA,AJ,BC,CR,D,F,IM-1,IM-2,MC1,MC2} and the references therein. However, the rotating trapped BECs in the attractive case behave different extremely from those in the well-understood repulsive case. Typically, the vortices are generally unstable in the attractive case (cf. \cite{CC,LC}), even though the vortices may form stable lattice configurations in the repulsive case, cf. \cite{A,F}.  Because of distinct mechanisms,  existing physical observations and numerical simulations show that rotating trapped BECs with attractive interactions present more complicated phenomena and structures, see \cite{BC,CC,D,F,LC}, only few of which have however been investigated analytically so far.

As derived rigorously in \cite{LNR,Lewin,NR,R,Roug} by a mean-field approximation, the ground state of two-dimensional attractive BECs in a rotating trap can be described equivalently by a complex constraint minimizer of the following Gross-Pitaevskii (GP) energy functional
\begin{equation}
F_a(u):=\int _{\R ^2} \big(|\nabla u|^2+V(x)|u|^2\big)dx-\frac{a}{2}\int _{\R ^2}|u|^4dx-
\Omega \int_{\R ^2}x^{\perp}\cdot (iu,\, \nabla u)dx,  \ \ u\in \h , \label{f}
\end{equation}
under the mass constraint
\begin{equation}\label{def:ea}
	e_F(a):=\inf _{\{u\in \h, \, \|u\|^2_2=1 \} } F_a(u),\ \, a>0,
\end{equation}
where $x^{\perp} =(-x_2,x_1)$ with $x=(x_1,x_2)\in \R^2$, $(iu,\, \nabla u)=i(u\nabla \bar u-\bar u\nabla u)/2$, and the complex space $\h$ is defined as
\begin{equation}\label{1.H}
\h :=  \Big \{u\in  H^1(\R ^2, \mathbb{C}):\ \int _{\R ^2}
V(x)|u|^2 dx<\infty\Big \}.
\end{equation}
The parameter $a>0$ in $e_F(a)$ denotes the product for the number $N$ of particles times the absolute value of the scattering length  $\nu$ in the two-body
interaction, and while $\Omega \ge 0$ describes the rotational velocity of the rotating trap $V(x)\ge 0$. We comment that one may impose $e_F(a)$ a different constraint $\int_{\mathbb{R}^2} |u(x)|^2dx=N>0$, but these two different forms can be   reduced equivalently to each other, see \cite{GLY}. In this paper we therefore focus on the form of $e_F(a)$ instead. It also deserves to remark that even though the mass-subcritical version of $e_F(a)$, where the nonlinear term $|u|^4$ of $F_a(u)$ is replaced by $|u|^p$ for $2<p<4$, was studied as early as in the pioneering work of Esteban-Lions \cite{EL}, the mass-critical constraint variational problem  $e_F(a)$ in the complex range was not addressed until recent years, see \cite{ANS,BC,BH,GLP,GLY,Lewin} and the references therein.

The non-rotational case $\Omega=0$ of $e_F(a)$ was studied recently in \cite{GLW,GS,GWZZ,LPY,Z} and the references therein, where the existence, uniqueness, symmetry breaking and other analytical properties of complex-valued  minimizers were investigated equivalently in view of the argument \cite[Theorem II.1]{CL}.
For a class of trapping potentials $V(x)$, above mentioned works give implicitly that
$e_F(a)$ with $\Omega=0$ admits complex-valued minimizers if and only if $a<a^*$, where $a^*=\|w\|^2_{L^2(\R^2)}$ and $w=w(|x|)>0$ is the unique (cf. \cite{K,W}) positive solution of the following nonlinear scalar field equation
\begin{equation}
\Delta u-u+u^3=0\  \mbox{  in } \  \R^2,\,\ u\in H^1(\R ^2,\R).  \label{Kwong}
\end{equation}

Starting from earlier works \cite{BC,Lewin}, the rotational case $\Omega>0$ of $e_F(a)$ was analyzed more recently. More precisely, considering special trapping potentials, such as typically $V(x)=|x|^2$, the existence and nonexistence, stability and some other properties of complex-valued minimizers for $e_F(a)$ with $\Omega>0$ were studied in \cite{ANS,BC,BH,Lewin}. Generally, if the trapping potential $0\leq V (x)\in L^{\infty}_{loc}(\R^2)$ satisfies
\begin{equation}\label{A:V}
\underline{\lim} _{|x|\to\infty }\frac{V(x)}{|x|^2}>0,
\end{equation}
then one can define as in \cite{GLY} the following critical rotational velocity $\Omega ^*:=\Omega ^*(V)$:
\begin{equation}
\Omega ^*:=\sup \Big\{\Omega >0:\ \  V(x)-\frac{\Omega ^2}{4}|x|^2  \to\infty \,\ \mbox{as}\,\
|x|\to\infty \Big\}.  \label{Omega}
\end{equation}
Note that if $V(x)$ satisfies the assumption (\ref{A:V}), then $\Omega^*\in (0,+\infty]$ exists, and $V_\Omega (x):=V(x)-\frac{\Omega^2}{4}|x|^2\ge 0$ holds in $\R^2$ for any $0\le\Omega <\Omega ^*$. Under the assumption (\ref{A:V}), we established in \cite[Theorem 1.1]{GLY} the following existence and non-existence of complex-valued minimizers:

\vskip 0.05truein
\noindent{\bf Theorem A.} {\em  (\cite[Theorem 1.1]{GLY}) Assume $V(x)\in L^\infty_{\rm loc}(\R^2)$ satisfies (\ref{A:V}) such that $\Omega^*\in (0,+\infty]$ in \eqref{Omega} exists. Then we have
\begin{enumerate}
		\item If $0\le\Omega <\Omega ^*$ and $0\leq a< a^*:=\|w\|^2_2$, then there exists at least
		one minimizer of $e_F(a)$.
		\item If $0\le\Omega <\Omega ^*$ and $a \ge a^*:=\|w\|^2_2$, then there
		is no minimizer of $e_F(a)$.
		\item If $\Omega >\Omega ^*$, then for any $a\ge 0$, there is no minimizer of $e_F(a)$.
\end{enumerate}}

The proof of Theorem A needs the following Gagliardo-Nirenberg inequality
\begin{equation}\label{11:GNineq}
	\inte |u(x)|^4 dx\le \frac 2 {\|w\|_2^{2}} \inte |\nabla u(x) |^2dx \inte |u(x)|^2dx ,\
\  u \in H^1(\R ^2, \R),
\end{equation}
where the identity is attained at $w$, and the following diamagnetic
inequality
\begin{equation}  |\nabla u|^2 -\Omega \, x^{\perp}\cdot (iu,\, \nabla u) =
	|(\nabla -i\mathcal{A } )u|^2 -\frac{\Omega ^2}{4}  |x|^2|u|^2\ge \big| \nabla |u|\big|^2
	-\frac{\Omega ^2}{4}  |x|^2|u|^2, \,\ u\in H^1(\R^2, \mathbb{C}),
	\label{1:2:1A}
\end{equation}
where $\mathcal{A }=\frac{\Omega}{2}x^{\perp}$, see \cite{Lieb,W} for more details. By the variational theory, if $e_F(a)$ admits a minimizer $u_a$, then $u_a$ is a ground state of
the following Euler-Lagrange  equation
\begin{equation}
-\Delta u_a+V(x)u_a+i\, \Omega \, (x^{\perp}\cdot \nabla u_a)=\mu u_a+a|u_a|^2u_a\quad \mbox{in}
\, \  \R^2, \ \ \inte |u_a|^2dx=1, \label{eqn}
\end{equation}
where $\mu = \mu (a, \Omega, u_a )\in \R$ is a suitable Lagrange multiplier. We comment that there exist many interesting progresses on normalized solutions of the elliptic problem (\ref{eqn}), see \cite{AS,BJ,BL,BW,Cao,EL,Gou,J,LPY,PP} and the references therein.


By employing energy estimates and elliptic PDE theory, it was proved in \cite{GLY,Lewin} that the minimizer $u_a$ of $e_F(a)$ concentrates at a global minimum point of $V_\Omega (x)$ as $a\nearrow a^*$, in the sense that
\begin{equation}
\|u_a\|_\infty\to \infty\ \ \mbox{and }\ \   \inte  V_\Omega (x)|u_a|^2dx\to V_\Omega (x_0):
=\inf _{x\in \R^2}
V_\Omega (x)\quad \mbox{as}\ \ a\nearrow a^*.
\label{intro:lower}
\end{equation}
Based on (\ref{intro:lower}), the $L^\infty$ uniform convergence of $u_a$  after rescaling and translation  was also obtained in \cite{GLY}.
By developing the method of inductive symmetry, we further proved in \cite[Theorem 1.3]{GLY} (see also \cite{LP}) the absence of vortices of minimizers $u_a$ for $e_F(a)$ as $a\nearrow a^*$ in the radially symmetric case $V(x)=|x|^2$, where the rotational velocity $\Omega$ presents essentially no effect on $u_a$ and $Im (u_ae^{i\delta_a})\equiv 0$ as $a\nearrow a^*$ for some constant phase $\delta_a\in [0, 2\pi)$. We should emphasize that the arguments of \cite[Theorem 1.3]{GLY} cannot however be extended to the non-radially symmetric case of $V(x)$, since one cannot expect generally $Im (u_ae^{i\delta_a})\equiv 0$ as $a\nearrow a^*$ for some constant phase $\delta_a\in [0, 2\pi)$, see also Remark \ref{rem2.1} below.


\subsection{Main Results}

As a continuation of \cite{GLY}, it is natural to ask whether the rotational velocity $\Omega$ has some visible effect on  the minimizers $u_a$ of $e_F(a)$ as $a\nearrow a^*$, provided that the trap $V(x)$ is non-radially symmetric. The main purpose of this paper is to address the above question.

Due to the physical relevance (cf. \cite{MC1,MC2,AA,CD,IM-1}), in this paper we mainly focus on the non-radially harmonic trap $V(x)\ge 0$ of the form
\begin{equation}\label{5:V}
 V(x)=x_1^2+\Lambda ^2x_2^2\ \ \mbox{for}\ \ x=(x_1, x_2)\in\R^2,\quad\mbox{where}\ \ 0<\Lambda <1,
\end{equation}
since the case  $\Lambda >1$ can be established similarly. As mentioned before, the case  $\Lambda =1$ of (\ref{5:V}) was already addressed in Theorem 1.3 and (1.18) of \cite{GLY}. We comment that Theorems \ref{thm2}  and \ref{thm1} below can be extended to the general case where the trap $V(x)$ is  homogeneous of degree $2$, see (\ref{1:V}) for its definition, and however for simplicity we do not pursue such extensions in this paper.

Under the assumption (\ref{5:V}), one can note from (\ref{Omega}) that $e_F(a)$ admits the critical velocity $\Omega ^*:=2\Lambda$. For any fixed $0<\Omega<\Omega ^*:=2\Lambda<2$, we further have
\begin{equation}\label{5:2}
\begin{split}
&\qquad (0, 0) \,\ \text{is the unique and non-degenerate critical point of}\\
 & H_\Omega (y):=\inte V_\Omega (x+y)w^2(x)dx,\ \ V_\Omega (x):=V(x)-\frac{\Omega^2}{4} |x|^2\ge 0,
\end{split}\end{equation}
and
\begin{equation}\label{5:2K}
V_\Omega(x)+\frac{\Omega ^2 }{4}|x|^2=V(x)=x_1^2+\Lambda ^2x_2^2.
\end{equation}
Under the assumption (\ref{5:V}), we thus note from \cite[Theorem 1.1]{GLP} that up to the constant phase, complex-valued minimizers of $e_F (a)$ must be unique as $a\nearrow a^*$.
Define
\begin{equation}
\label{5:3}
\lambda _0  =\Big(   \inte V(x)w^2dx  \Big)
^{\frac{1}{4}}>0.
\end{equation}
For convenience, we denote  $\varphi _i(x)\in C^2(\R^2)\cap L^\infty(\R^2)
$ to be the unique solution of
\begin{equation}\label{thm2:1}
 \nabla \varphi _i(0)=0,\quad \big(-\Delta+1-3w^2\big)\varphi _i(x)=\tilde{f}_i(x)
 \ \ \mbox{in}\,\ \R^2,\ \ i=1,\, 2,
\end{equation}
where $\tilde{f}_i(x)$ satisfies
\begin{equation}\label{thm2:11}\arraycolsep=1.5pt
\tilde{f}_i(x)=\left\{\begin{array}{lll}
&- \displaystyle\frac{\lam _0^4}{a^*}w^3(x)-\big(x_1^2+\Lambda ^2x_2^2\big) w(x),\quad \quad \ &\mbox{if}\ \ i=1;\\[3mm]
 &- (1-\Lambda ^2)x_1^2 w(x), \,\ &\mbox{if}\ \ i=2.
\end{array}\right.\end{equation}
Here the uniqueness of $\varphi_{i}(x)$ follows from   \cite[Theorem 1.2]{Frank} and \cite[Lemma 4.1]{Wei96}  together with the fact $\nabla \varphi_{i}(0)=0$ for $i=1, 2$.
We also denote $\varphi_{I}(x)\in C^2(\R^2)\cap L^\infty(\R^2)$ to be the unique solution of
\begin{equation}\label{thm2:2}
 \big(-\Delta+1-w^2\big) \varphi_{I}(x)=-  \big(x^\perp \cdot \nabla \varphi_2\big)
 \ \ \mbox{in}\,\ \R^2,\ \,  \inte \varphi_{I}wdx=0,
\end{equation}
where $\varphi_2$ satisfies (\ref{thm2:1}). Here the uniqueness of $\varphi_{I}(x)$ follows from \cite[Lemma 4.1]{Wei96} and the fact $\int_{R^2} \varphi_{I}wdx=0$. Applying above notations, {\em the main result of this paper} is concerned with the following asymptotic expansion as $a\nearrow a^*$, which displays the visible effect of the rotational velocity $\Omega$:

\begin{thm}\label{thm2}
Suppose that the non-radially harmonic trap $V(x)$ satisfies (\ref{5:V}) for some $0<\Lambda <1$. Then for any fixed $0<\Omega <\Omega ^*:=2\Lambda$, the complex-valued minimizer $u_a$ of $e_F(a)$ satisfies
\begin{equation}\arraycolsep=1.5pt\begin{array}{lll}
v_a:=\ep  \sqrt{a^*}\, u_{a}\big( \ep  x +x_a \big)e^{-i\, \big(\frac{\Omega}{2 }\ep x\cdot x_a^\perp-\theta_a\big)}& =&
w+\ep^4\big\{\varphi _1+C_\Lambda(w+x\cdot \nabla w)\big\}\\[2mm]
&&+i\,\Omega\ep^6   \varphi_{I}+o(\ep^6) \ \ \mbox{in}\ \ L^\infty(\R^2, \mathbb{C})
\end{array}
\label{thm2:3}
\end{equation}
as $a\nearrow a^*$, where $\ep := \frac{(a^*-a)^{\frac{1}{4}}}{\lambda _0}>0$, $x_a$ is the unique global maximum point of $|u_a|$ and satisfies
\begin{equation}
\big|   x_a \big|= o\big( \ep ^5\big)\ \ \mbox{as} \ \ a\nearrow a^*,
\label{thm2:4}
\end{equation}
$\theta_a \in [0,2\pi)$ is a suitable constant phase, and the constant $C_\Lambda\not =0$ is independent of $\Omega$ and satisfies
\begin{equation}\label{thm2:5}
  C_\Lambda=\frac{1 }{2\lam_0 ^4}\Big[\displaystyle\inte (3w^2-1)\varphi^2_1- \displaystyle 4\inte \big(x_1^2+\Lambda ^2x_2^2\big)w \varphi_1\Big]\not =0.
\end{equation}
Here $\varphi _1$ and $\varphi _I$ are uniquely given by (\ref{thm2:1})-(\ref{thm2:2}).
\end{thm}

\begin{rem}
\label{rem4.1}
(1). Under the assumption (\ref{5:V}), the interesting novelty of Theorem \ref{thm2} lies in the fact that  the rotational velocity $\Omega$ affects visibly the minimizer $u_a$ starting from  the third term, which is imaginary, of $v_a$ defined in (\ref{thm2:3}). This is however different from the radial case $V(x)=|x|^2$ addressed earlier in \cite[Theorem 1.3]{GLY}, where $\Omega>0$ has essentially no effect on $u_a$.



(2). Specially, if $V(x)=|x|^2$, then the argument of Theorem \ref{thm2} can yield that (\ref{thm2:3}) holds for $\varphi_{I}\equiv 0$, and while the term $\varphi _1+C_\Lambda(w+x\cdot \nabla w)$ is radially symmetric and  independent of $\Omega>0$. Following these, we guess that one may further obtain the non-existence of vortices for $u_a$ in $\R^2$ as $a\nearrow a^*$, which approach is different slightly from that of \cite{GLY}. We leave it to the interested reader.

(3). We expect that, instead of (\ref{thm2:4}), the unique global maximum point $x_a$ of $|u_a|$ satisfies $x_a\equiv 0$ as $a\nearrow a^*$. But our analysis does not give this conclusion, for which one needs to derive further refined estimates of $u_a$  as $a\nearrow a^*$.
\end{rem}

We next follow three steps to explain the general strategy of proving Theorem \ref{thm2}. As the first step, we proved in \cite[Theorem 1.2]{GLY} that the complex-valued minimizer $u_a$ of $e_F(a)$ satisfies the equation (\ref{eqn}) and
\begin{equation}
v_a(x):=\ep  \sqrt{a^*}\, u_{a}\Big( \ep  x +x_a \Big)e^{-i\, \big(\frac{\Omega}{2 }\ep x\cdot x_a^\perp-\theta_a\big)}:=\big[R_a(x)+w(x)\big]+iI_a(x) \to   w(x)
\label{0:v:1}
\end{equation}
uniformly in $L^\infty (\R^2, \mathbb{C})$ as $a\nearrow a^*$, where and below $\ep  >0$ and $x_a\in \R^2$ are as in Theorem \ref{thm2}, and $\theta_a \in [0, 2\pi )$ is chosen suitably. Under the assumption  (\ref{5:V}), recall from \cite[Section 3]{GLY} that the Lagrange multiplier $\mu_{a}$ of (\ref{eqn}) and $x_a$ satisfy
\begin{equation}
 \mu_{a}\ep ^2\to -1 \ \ \mbox{and}\ \ \frac{x_a}{\ep }\to 0\ \ \mbox{as} \ \ a\nearrow a^*.
\label{0:v:2}
\end{equation}
Note from (\ref{0:v:1}) that $R_a$ and $I_a$ satisfy
 \begin{equation}\label{0:2.2c}
  R_a(x)\to 0\  \    \mbox{and}\ \ I_a(x)\to 0\,\ \mbox{uniformly in}\,\ \R^2\ \ \mbox{as} \ \ a\nearrow a^*.
 \end{equation}

The second step of proving Theorem \ref{thm2} is to establish the  leading terms of $R_a(x)+iI_a(x)$, in terms of $\ep $ and $\mu_{a}\ep ^2+1$, by following (\ref{0:v:1}) and (\ref{0:v:2}). However, we remark that it seems difficult to reach this aim by investigating directly the Euler-Lagrange equation (\ref{eqn}). To overcome this difficulty, motivated by (\ref{eqn}) and (\ref{0:v:1}) we shall consider the following coupled system of  $R_a(x)$ and $I_a(x)$:
\begin{equation}\label{0:2.21a}
\arraycolsep=1.5pt\left\{\begin{array}{lll}
\big(\mathcal{L}_a -w^2-w (R_a+w)\big) R_a&=&\ep ^2 \Omega\, (x^\bot\cdot\nabla I_a)+\tilde{F}_a(x)
\ \ \mbox{in}\,\  \R^2,\ \ \nabla R_a(0)\not \equiv 0,\\ [3mm]
\qquad\qquad\qquad\qquad\quad\, \  \mathcal{L}_a I_a&=&-\ep ^2 \Omega\, (x^\bot\cdot\nabla R_a)
\qquad\quad \ \mbox{in}\,\  \R^2,\quad \displaystyle\inte wI_adx\equiv 0,
\end{array}\right.
\end{equation}
where  $\tilde{F}_a(x)$ is an inhomogeneous term containing $ \mu_{a}\ep  ^2+1 $, and the operator $
\mathcal{L}_a$ is as in (\ref{2.20}).

Even though the rest part of the second step is motivated by \cite[Theorem 1.4]{GLW} and (\ref{0:2.2c}), unfortunately, there appear extra difficulties. Actually, we first note from (\ref{0:2.21a}) that we have $R_a(0)\not \equiv 0$ in $\R^2 $, which leads to a new difficulty in studying the expansion of $R_a$, see Lemma \ref{lem2.2}. Moreover, the coupled system (\ref{0:2.21a}) contains the coupled rotating terms $(x^\bot\cdot\nabla I_a)$ and $(x^\bot\cdot\nabla R_a)$. To overcome these extra difficulties, we shall proceed with the refined estimates of $R_a$ and $I_a$, and make full use of the non-degenerancy of $w$ as well. We shall finally obtain in Lemmas \ref{lem2.3} and \ref{lem3.1} the leading terms of $I_a$ and $R_a$, respectively, in terms of $\ep $ and $\mu_{a}\ep ^2+1 $.

The third step of proving Theorem \ref{thm2} is to address the refined estimate of  $\mu_{a}\ep ^2+1$ in terms of $\ep >0$. We shall achieve it by taking full advantage of the mass constraint $\inte |u_a|^2dx=1$, as well as analytical results of the second step.

Consider the non-radially harmonic trap $V(x)$ satisfying (\ref{5:V}). It then follows from \cite[Theorem 1.2]{GLY} that the minimizers do not have any vortex near the origin as $a\nearrow a^*$.
As a byproduct of Theorem  \ref{thm2}, we shall derive the following nonexistence of vortices under the assumption (\ref{5:V}).

\begin{thm}\label{thm1}
Suppose  the non-radially harmonic trap $V(x)$ satisfies (\ref{5:V}),
and let $u_a$ be a complex-valued minimizer of $e_F(a)$, where $0<\Omega <\Omega ^*:=2\Lambda$ is fixed.
Then there exists a constant $C>0$, independent of $0<a<a^*$, such that
\begin{equation}
|u_a(x)|>0\,\  \mbox{in the region}\,\ R(a):=\big\{x\in\R^2:\, |x|\le C(a^*-a)^{-\frac{1}{12}} \big\}\,\  \mbox{as}\,\ a\nearrow a^*,
\label{thm1:2}
\end{equation}
i.e., $u_a$ does not admit any vortex in the region $R(a)$  as $a\nearrow a^*$.
\end{thm}

\begin{rem}
\label{rem1:B}
(1). The nonexistence of vortices for rotating trapped BECs was studied earlier by jacobian estimates, vortex ball constructions, the inductive symmetry, and some other methods,  see
\cite{AJ,AN,IM-1,IM-2,AA,CD,CR,GLY,LP,SSbook} and the references therein. However, as far as we know,  above mentioned methods depend heavily on the radial symmetry of $V(x)$.

(2). Theorem \ref{thm1} proves the nonexistence of vortices  in a very large region $R(a)$ as $a\nearrow a^*$, including the places where $|u_a|$ is already very small.
As remarked before, the argument of Theorem \ref{thm1}  holds actually for {\em the non-radially general case} where the trap $V(x)$ is  homogeneous of degree $2$, see (\ref{1:V}) for its definition. This is the main contribution of Theorem \ref{thm1}, since the existing methods of \cite[Theorem 1.1]{AJ} and \cite[Theorem 1.3]{GLY} and the references therein are not applicable to the non-radially symmetric case of $V(x)$.

%
\end{rem}

To prove Theorem \ref{thm1} with the non-radially symmetric trap $V(x)$, we shall make full use of  Theorem \ref{thm2} to establish the following optimal estimate
\begin{equation}\label{1:OK:9}
    \big|v_{a}(x)-w(x)\big|
    \leq C_1 \ep ^{4}|x|^{\frac{5}{2}}e^{-\sqrt{1-C_2\ep ^{4}}|x|}\ \ \hbox{uniformly in}\ \ \R^2 \ \ \mbox{as} \ \ a\nearrow  a^*,
\end{equation}
where $v_{a}(x)$ is as in (\ref{thm2:3}), and the constants $C_1>0$ and $C_2>0$ are independent of $0<a<a^*$. Applying (\ref{1:OK:9}), Theorem \ref{thm1} then follows from  the fact that $|v_{a}|\geq |w|-|w-v_{a}|$ holds in $\R^2$, due to the exact exponential decay (\ref{2:exp}) of $w$.
We note from above that the idea of proving Theorem 1.2 is different from those of \cite[Theorem 1.3]{GLY} and \cite[Theorem 1.3]{GLY}, since the latter ones proved the nonexistence of vortices in the following approach: $v_{a}(x)$ is essentially a real-valued minimizer of $e_F(a)$ with $\Omega =0$, and $|v_{a}(x)|>0 $ then holds in a standard way, see \cite{GS,GWZZ}.

%


This paper is organized as follows: In Section 2 we study refined estimates of minimizers $u_a$ for $e_F(a)$ as $a\nearrow a^*$. Applying analytical results of the previous section, in Section 3 we first analyze the refined estimate of  $\mu_{a}\ep ^2+1$ in terms of $\ep >0$, based on which we then finish the proofs of   Theorems \ref{thm2} and \ref{thm1} in Subsection 3.1. In Appendix A we shall address the proof of the  claims (\ref{5:8}) and (\ref{5:2:beta}) used in Section 3.

\section{Refined Estimates of Minimizers}
Suppose the rotating speed $\Omega\in (0, \Omega^*)$ is fixed, the purpose of this section is to address the refined estimates of complex-valued minimizers $u_a$ for $e_F(a)$  as $a\nearrow a^*$. To clarify the general idea, in this section we focus on the homogeneous trapping potential $V(x)$, in the sense that a function $h(x):\R^2\longmapsto \R$ is called homogeneous of degree $p>0$ (about the origin), if
\begin{equation}\label{1:V}
h(tx)=t^ph(x)\,\ \hbox{for any $t\in\R^+$ and $x\in \R^2$.}
\end{equation}
We note from above that if  $0\le V(x)\in C^2(\R^2)$  is homogeneous of degree $2$ and satisfies  $\lim_{|x|\to\infty} V(x) = \infty$, such as (\ref{5:V}), then $x=0$ is the unique minimum point of $V(x)$, and
\begin{equation}\label{1:VO}
0\le V_\Omega (x):=V(x)-\frac{\Omega^2}{4} |x|^2\in C^2(\R^2)\ \, \mbox{is also homogeneous of degree $2$ }
\end{equation}
for any fixed $0< \Omega<\Omega ^*$, where  $0<\Omega^*<\infty$ is defined as in (\ref{Omega}).

Recall from \cite{W} that the unique positive radial solution $w$ of \eqref{Kwong} is an optimizer of the following
Gagliardo-Nirenberg inequality
\begin{equation}\label{GNineq}
\inte |u(x)|^4 dx\le \frac 2 {a^*} \inte |\nabla u(x) |^2dx \inte |u(x)|^2dx ,\
\  u \in H^1(\R ^2, \R).
\end{equation}
Note also from \cite[Lemma 8.1.2]{C} and \cite[Proposition 4.1]{GNN} that $w=w(|x|)>0$ satisfies
\begin{equation}\label{2:id}
\inte |\nabla w |^2dx  =\inte w ^2dx=\frac{1}{2}\inte w ^4dx,
\end{equation}
and
\begin{equation}
w(x) \, , \ |\nabla w(x)| = O(|x|^{-\frac{1}{2}}e^{-|x|}) \quad
\text{as \ $|x|\to \infty$.}  \label{2:exp}
\end{equation}
Define for $0<a<a^*$,
\begin{equation}\label{2.e}
\eps_a:= \frac{(a^*-a)^{\frac{1}{4}}}{\lambda}>0,
\quad
\lambda   =\Big[   \inte \Big(V_\Omega(x+y_0)+\frac{\Omega ^2 }{4}|x|^2\Big)w^2(x)dx  \Big]
^{\frac{1}{4}}>0,
\end{equation}
where $V_\Omega (x)\ge 0$ is defined by (\ref{1:VO}), and $y_0\in\R^2$ denotes a unique global minimum point of $ H_\Omega (y):=\inte V_\Omega (x+y)w^2(x)dx$.
Define
\begin{equation}
\label{2:v}
v_a(x):=\varepsilon_a \sqrt{a^*}\, u_a\big(\varepsilon_a x+x_a\big)e^{-i (\frac{\varepsilon _a \Omega}{2} x\cdot x_{a}^{\bot}-\theta_a)}=\widetilde{R}_a(x)+iI_a(x),
\end{equation}
where $x_a$ is a global maximal point of $|u_a(x)|$, $\widetilde{R}_a(x)$ and $I_a(x)$ denote the real and imaginary parts of $v_{a}(x)$, respectively, and the constant phase $\theta_a\in [0,2\pi)$ is chosen such that
\begin{equation}\label{2.theta}
\big\| v_{a}-w\big\|_{L^2(\R^2)}=\min_{\theta\in [0,2\pi)}\big\|e^{i\theta} \tilde{v}_{a}-w\big\|_{L^2(\R^2)},
\end{equation}
where $\tilde{v}_{a}:=\varepsilon_a \sqrt{a^*}\, u_a\big(\varepsilon_a x+x_a\big)e^{-i \frac{\varepsilon _a \Omega}{2} x\cdot x_{a}^{\bot}}$.
The above property gives the following orthogonality condition on $I_a(x)$:
\begin{equation}\label{2:I}
\inte w(x)I_a(x)dx=0.
\end{equation}
Using above notations, {we proved in \cite{GLY} the following  $L^\infty-$uniform convergence as $a\nearrow a^*$}.

\begin{prop}\label{prop2.1}
(\cite[Theorem 1.2]{GLY}) Assume $0\le V(x)\in C^2(\R^2)$ satisfying $\lim\limits_{|x|\to\infty} V(x) = +\infty$ is homogeneous of degree $2$, and let $0$ be a unique global minimum point of $ H_\Omega (y):=\inte V_\Omega (x+y)w^2(x)dx$. For any fixed $0<\Omega <\Omega ^*$, where $\Omega^*>0$ is defined as in (\ref{Omega}), suppose  $u_a$ is a minimizer of $e_F(a )$. Then we have
\begin{equation}
v_a(x):=\eps_a \sqrt{a^*}\, u_{a}\Big( \eps_a x +x_a \Big)e^{-i\, \big(\frac{\Omega}{2 }\eps_ax\cdot x_a^\perp-\theta_a\big)} \to   w(x)
\label{2:v:1}
\end{equation}
uniformly in $L^\infty (\R^2, \mathbb{C})$ as $a\nearrow a^*$, where $\eps_a>0$ is as in (\ref{2.e}), $\theta_a \in [0,2\pi)$ is  chosen such that (\ref{2.theta}) holds true, and the unique global maximal point  $x_a\in\R^2$ of $|u_{a}|$ satisfies
\begin{equation}
\lim _{a\nearrow a^*}\frac{x_a}{\eps_a}=0.
\label{2:v:2}
\end{equation}
\end{prop}

By the variational theory, the minimizer $u_a$ of $e_F(a )$ solves the following Euler-Lagrange equation:
\begin{equation}\label{2.1}
-\Delta u_a+V(x)u_a+i\,\Omega (x^{\perp}\cdot \nabla u_a)=\mu_au_a+a|u_a|^2u_a\,\ \,\hbox{in}\,\ \R^2,				 
\end{equation}
where $\mu_a:=\mu_a(u_a)\in \R$ is a suitable Lagrange multiplier satisfying
\begin{equation}\label{2.2}
\mu_a=e_F(a)-\frac{a}{2}\inte |u_a|^4dx.
\end{equation}
Note from (\ref{2.1}) that the function $v_a$ satisfies the following elliptic equation
\begin{equation}\label{2:va}
\begin{split}
&-\Delta v_a+i\,\eps_a^2\,\Omega \big(x^{\perp}\cdot\nabla v_a\big)+\Big[\frac{\eps_a^4\Omega^2|x|^2}{4}+\eps_a^2V_{\Omega}(\eps_ax
+x_{a})\Big]v_a\\
=&\mu_{a}\eps_a^2v_a+\frac{a}{a^*}|v_a|^2v_a\quad\hbox{in }\,\,\R^2,
\end{split}
\end{equation}
where $V_\Omega (x)\ge 0$ is as in (\ref{1:VO}).
The analysis of \cite[Theorem 1.2]{GLY} gives that the above Lagrange multiplier $\mu_a$ satisfies
\begin{equation}
\lim _{a\nearrow a^*} \mu_a\eps_a^2=-1.
\label{2:v:3}
\end{equation}
In order to obtain the refined estimate of $u_a$ as $a\nearrow a^*$, one can note from (\ref{2:v:1}) and (\ref{2:va}) that a more refined estimate than (\ref{2:v:3}) is needed for $\mu_a$ as $a\nearrow a^*$, which is one of the main difficulties in this paper. Towards this purpose, it however seems difficult to handle directly with the single equation (\ref{2:va}), instead of which we shall consider the coupled system of $Re(v_a)$ and $Im(v_a)$ in the coming subsection.

\subsection{Refined estimates of $v_a$ as $a\nearrow a^*$}

Based on Proposition \ref{prop2.1}, in this subsection we shall derive refined estimates of $v_a$ defined in (\ref{2:v}) as $a\nearrow a^*$. Towards this purpose, we introduce the following linear operator
\begin{equation}\label{2:2:1}
\mathcal{L}:=-\Delta+1-w^2\quad\hbox{in}\,\ \R^2.
\end{equation}
It then obtains from \cite[Theorem 11.8]{Lieb} and \cite[Corollary 11.9]{Lieb}  that
\begin{equation}\label{2:2:2}
ker\mathcal{L}=\{w\}\quad\hbox{and}\quad \langle\mathcal{L}v,v\rangle\geq 0 \ \ \mbox{for all} \ \ v\in L^2(\R^2).
\end{equation}
We also define the linearized operator $\widetilde{\mathcal{L}}$ of (\ref{Kwong}) around $w>0$ by
\begin{equation}\label{2:2:3}
\widetilde{\mathcal{L}}:=-\Delta+1-3w^2\quad\hbox{in}\,\ \R^2.
\end{equation}
It follows from \cite{K,NT,Wei96} that
\begin{equation}\label{2:2:4}
ker\widetilde{\mathcal{L}}=\Big\{\frac{\partial w}{\partial x_1},\ \frac{\partial w}{\partial x_2}\Big\}.
\end{equation}
For convenience, we denote for $0<a<a^*$,
\begin{equation}\label{2:beta}
\eps _a:= \frac{(a^*-a)^{\frac{1}{4}}}{\lambda}>0,\,\ \alp _a:=(\lambda \eps_a) ^4=a^*-a >0, \ \ \mbox{and} \ \  \beta _a:=1+ \mu _a\eps _a^2 ,
\end{equation}
where $\lam>0$ is defined by (\ref{2.e}) with $y_0=(0, 0)$ and the Lagrange multiplier $\mu _a$ satisfies (\ref{2:v:3}). We then have
\[\alp _a \to 0 \ \ \mbox{and} \ \ \beta _a\to 0   \,\ \mbox{as} \,\ a\nearrow a^*. \]

Following (\ref{2:v}), we now rewrite $v_a$ as
\begin{equation}\label{2:R}
v_a(x):=\widetilde{R}_a(x)+iI_a(x)=\big[R_a(x)+w(x)\big]+iI_a(x),
\end{equation}
so that
\[
R_a(x)\to 0\,\ \hbox{and}\,\ I_a(x)\to 0\,\ \hbox{uniformly in}\,\ L^\infty (\R^2, \R) \ \ \hbox{as} \,\ a\nearrow a^*,
\]
due to Proposition \ref{prop2.1}. Since $\nabla |v_a(0)|\equiv 0$ holds for all $0<a<a^*$, we derive from (\ref{2:R}) that
\begin{equation}\label{4:3RT}
\nabla R_a(0)=-\frac{I_a(0)\nabla I_a(0)}{ w(0) +R_a(0)}\to 0  \  \ \mbox{as}\ \ a\nearrow a^*.
\end{equation}
For simplicity, we denote the operator $\mathcal{L}_a$ by
\begin{equation}\label{2.20}
\mathcal{L}_a:=-\Delta+\Big(\frac{\eps_a^4\Omega^2}{4}|x|^2+\eps_a^2V_\Omega (\eps_a x+x_a)-\mu_a \e^2 -\frac{a}{a^*}|v_a|^2\Big)\quad\hbox{in }\,\,\R^2.
\end{equation}
It then follows from (\ref{2:I}),  (\ref{2:va}) and (\ref{2:R}) that $ I_a$ satisfies
\begin{equation}\label{2.21}
\mathcal{L}_a I_a=-\eps^2_a\Omega\, (x^\bot\cdot\nabla R_a)
\ \ \mbox{in}\,\  \R^2,\quad \inte wI_adx\equiv 0,
\end{equation}
and while $R_a$ satisfies (\ref{4:3RT}) and
\begin{equation}\label{2.21a}
\widetilde{\mathcal{L}}_aR_a:=\Big[\mathcal{L}_a -w^2-w \widetilde{R}_a \Big]R_a=F_a(x)
\ \ \mbox{in}\,\  \R^2,
\end{equation}
where $F_a(x)$ is defined by
\begin{equation}\label{2.21b}\begin{split}
F_a(x):=&\eps^2_a\,\Omega\, (x^\bot\cdot\nabla I_a)-\Big[\frac{\eps_a^4\Omega^2}{4}|x|^2+\eps_a^2V_\Omega (\eps_a x+x_a)-\beta _a-\frac{a}{a^*}I_a^2+\frac{\alp_a}{a^*}\widetilde{R}^2_a\Big]w\\
=&-\eps^4_a\Big[\frac{\Omega^2}{4}|x|^2+V_\Omega \Big(x+\frac{x_a}{\eps_a}\Big)\Big]w+\beta _aw-\frac{\alp_a}{a^*}\widetilde{R}^2_a w\\
& +\eps^4_a\,\Omega\, \Big(x^\bot\cdot\frac{\nabla I_a}{\eps^2_a}\Big)+\frac{a}{a^*}I_a^2w.
\end{split}\end{equation}

Similar to \cite[Lemma 4.2 (1)]{GLY}, where $\eps_a>0$ is defined in a slightly different way, one can deduce from (\ref{2:va}) that there exists a constant $C>0$, independent of $0<a<a^*$, such that as $a\nearrow a^*$,
\begin{equation}\label{2.22}
  |\nabla  \widetilde{R}_a (x)|, \ |\nabla  I_a (x)|\leq Ce^{-\frac{2}{3}|x|}\quad \hbox{uniformly in}\, \ \R^2.
\end{equation}
Applying (\ref{4:3RT}) and (\ref{2.22}), the argument of \cite[Lemma 4.3]{GLY} then yields from (\ref{2.21}) that
there exists a constant $C>0$, independent of $0<a<a^*$, such that as $a\nearrow a^*$,
\begin{equation}\label{2.23}
|\nabla  I_a (x)|, \ |I_a (x)|\leq C \eps_a^2e^{-\frac{1}{8}|x|}\quad \hbox{uniformly in}\, \ \R^2.
\end{equation}
Setting
$$\psi_{2}(|x|):=-\frac{1}{2}\big(w+x\cdot \nabla w\big),$$
which is independent of $0<\Omega<\Omega ^*$, one can check that $\psi_{2}$ is a unique solution of
\begin{equation}\label{2.23A}
\nabla \psi_{2}(0)=0,\quad \widetilde{\mathcal{L}}\psi_{2}(|x|)=w(x)
 \,\ \mbox{in}\,\ \R^2,
\end{equation}
where the operator $\widetilde{\mathcal{L}}$ is defined by (\ref{2:2:3}). Having above estimates, we next establish the following ``rough"  limit profiles in terms of $\eps_a$ and $\beta_a$.

\begin{lem}\label{lem2.2}
Under the assumptions of Proposition \ref{prop2.1}, we have
\begin{enumerate}
 \item There exists a  constant $C>0$, independent of $0<a<a^*$, such that the imaginary part $I_a$  of (\ref{2:R}) satisfies
\begin{equation}\label{2.3:3A}
|\nabla  I_a(x)|, \ |I_a (x)|\leq C \eps_a^2(\e ^4+\gamma_a)e^{-\frac{1}{16} |x|}\quad \hbox{uniformly in}\, \ \R^2\,\ \mbox{as} \,\ a\nearrow a^*,
\end{equation}
where $\gamma_a>0$  satisfies $\gamma_a=o(|\beta_a|)$ as $a\nearrow a^*$.

\item  The real part $R_a$ of (\ref{2:R}) satisfies
\begin{equation}\label{2N:a3}
R_a(x):=\e ^4\psi_{1}(x)+\beta _a\psi_{2}(|x|) +o(\e ^4+|\beta _a|)\quad \hbox{in}\, \ \R^2 \,\ \mbox{as} \,\ a\nearrow a^*,
\end{equation}
where $\psi_{2}(|x|):=-\frac{1}{2}\big(w+x\cdot \nabla w\big)$ is radially symmetric, and $\psi_{1}(x)\in C^2(\R^2)\cap L^\infty(\R^2)$ solves uniquely
\begin{equation}\label{lem2.1:2}
 \nabla \psi_{1}(0)=0,\quad \widetilde{\mathcal{L}}\psi_{1}(x)=- \frac{\lam ^4}{a^*}w^3(x)-V(x)w(x)
 \,\ \mbox{in}\,\ \R^2.
\end{equation}
Here $\psi_{1}$ and $\psi_{2}$ are independent of $0<\Omega<\Omega ^*$.
\end{enumerate}
\end{lem}

\noindent{\bf Proof.} Denote
\begin{equation}\label{2.24R}
\mathcal{R}_a(x):=R_a(x)-\e ^4\psi_{1}(x)-\beta _a\psi_{2}(|x|),
\end{equation}
where $\psi_{2}(|x|):=-\frac{1}{2}\big(w+x\cdot \nabla w\big)$ is radially symmetric, and $\psi_{1}(x)\in C^2(\R^2)\cap L^\infty(\R^2)$ is a
solution of (\ref{lem2.1:2}). Note that $\psi_{1}$ and $\psi_{2}$ are independent of $0<\Omega<\Omega ^*$, and the uniqueness of $\psi_{1}(x)$ follows from
$\nabla \psi_{1}(0)=0$ and the property (\ref{2:2:4}) (see also \cite[Lemma 4.1]{Wei96} and \cite[Theorem 1.2]{Frank}). Moreover, by the comparison principle, see (\ref{2:D12}) below for the detailed argument, we derive from (\ref{2:exp}) and (\ref{lem2.1:2}) that
\begin{equation}\label{2NN:1}
\big|\psi_{1}(x)\big|,\  \big|\psi_{2}(|x|)\big|\le Ce^{-\delta|x|}\ \ \mbox{in}\ \ \R^2, \ \ \mbox{where}\ \ \frac{4}{5}<\delta <1.
\end{equation}
It then yields from (\ref{4:3RT}) and (\ref{lem2.1:2}) that $\mathcal{R}_a(x)$ satisfies
\begin{equation}\label{2.24a}
\nabla \mathcal{R}_a(0) \to 0   \, \ \mbox{as}\ \ a\nearrow a^*.
\end{equation}
We next carry out the proof by the following three steps:

{\em Step 1.}  Note from (\ref{2.21a}) and (\ref{2.24R}) that
\begin{equation}\label{2.24b}
\widetilde{\mathcal{L}}_a\mathcal{R}_a=\mathcal{F}_a(x)-\big(\widetilde{\mathcal{L}}_a-\widetilde{\mathcal{L}}\big)\big[\e ^4\psi_{1}(x)+\beta _a\psi_{2}(|x|)\big]:=N_a(x)\quad \hbox{in}\, \ \R^2.
\end{equation}
We deduce from (\ref{2.21b}), (\ref{2.23}) and (\ref{2NN:1}) that  the term $\mathcal{F}_a(x)$ satisfies
\[
 |\mathcal{F}_a(x)|:=\Big|F_a(x)-\beta_aw+\e^4\Big(\frac{\lam ^4}{a^*}w^3 +V(x)w \Big)\Big| \le C \e ^4  e^{-\frac{1}{10}|x|}\ \ \mbox{uniformly in}\ \ \R^2
\]
as $a\nearrow a^*,$
and
\[
\frac{\big|\big(\widetilde{\mathcal{L}}_a-\widetilde{\mathcal{L}}\big)\big[\e ^4\psi_{1}(x)+\beta _a\psi_{2}(|x|)\big]\big|}{\e ^4+ |\beta_a|}\le C\delta _a  e^{-\frac{1}{10}|x|}\ \ \mbox{uniformly in}\ \ \R^2\,\ \hbox{as}\, \ a\nearrow a^*  ,
\]
where  $\delta_a>0$  satisfies $\delta_a=o(1)$ as $a\nearrow a^*$.
We thus derive from above that
\begin{equation}\label{2NN:3}
\frac{|N_a(x)|}{\e ^4+ |\beta_a|}\le C  e^{-\frac{1}{10}|x|}\ \ \mbox{uniformly in}\ \ \R^2\,\ \hbox{as}\, \ a\nearrow a^*.
\end{equation}

We claim that there exists a constant $C>0$, independent of $0<a<a^*$, such that
\begin{equation}\label{2.24c}
\big| \mathcal{R}_a(x)\big| \leq C (\eps_a^4+|\beta_a|)  \quad \hbox{uniformly in}\, \  \R^2\,\ \hbox{as}\, \ a\nearrow a^*.
\end{equation}
Instead, assume that the above claim (\ref{2.24c}) is false, i.e., $\lim_{a\nearrow a^*}\frac{\|\mathcal{R}_a\|_{L^\infty(\R^2)}}{\eps_a^4+|\beta_a|}=\infty.$ Denote $\mathcal{\bar R}_a:=\frac{\mathcal{R}_a}{\|\mathcal{R}_a\|_{L^\infty(\R^2)}}$, so that $\|\mathcal{\bar R}_a\|_{L^\infty(\R^2)}=1$.
Applying (\ref{4:3RT}) and (\ref{2.23}), we derive from (\ref{2.24R}) that
\[
|\nabla \mathcal{\bar R}_a(0)|=\frac{|\nabla R_a(0)|}{\|\mathcal{R}_a\|_{L^\infty(\R^2)}}=\frac{1}{| w(0) +R_a(0)|}\frac{|I_a(0)||\nabla I_a(0)|}{\|\mathcal{R}_a\|_{L^\infty(\R^2)}}
\le \frac{C(\eps_a^4+|\beta_a|)}{\|\mathcal{R}_a\|_{L^\infty(\R^2)}}\, \ \mbox{as}\ \ a\nearrow a^*,\]
which then implies that
\begin{equation}\label{2.24dd}
\nabla \mathcal{\bar R}_a(0) \to 0   \ \   \mbox{as}\ \ a\nearrow a^*.
\end{equation}
We also deduce from (\ref{2.24b}) that
\begin{equation}\label{2.24d}
\widetilde{\mathcal{L}}_a\mathcal{\bar R}_a=\frac{N_a(x)}{\|\mathcal{R}_a\|_{L^\infty(\R^2)}}  \quad \hbox{in}\, \ \R^2.
\end{equation}
Note from (\ref{2NN:3}) that
\begin{equation}\label{2.24e}
\arraycolsep=1.5pt\begin{array}{lll}
\displaystyle \frac{|N_a(x)|}{\|\mathcal{R}_a\|_{L^\infty(\R^2)}} &:=&
\displaystyle\frac{\eps_a^4+|\beta_a|}{\|\mathcal{R}_a\|_{L^\infty(\R^2)}} \,\frac{\big|N_a(x)\big|}{\eps_a^4+|\beta_a|}\\[2mm]
&\le &\displaystyle\frac{\eps_a^4+|\beta_a|}{\|\mathcal{R}_a\|_{L^\infty(\R^2)}}C   e^{-\frac{1}{10}|x|}\quad \hbox{uniformly in}\, \  \R^2\,\ \hbox{as}\, \ a\nearrow a^*.
 \end{array}
\end{equation}
Suppose $y_a$ is a global maximum point of $|\mathcal{\bar R}_a(x)|$, so that $|\mathcal{\bar R}_a(y_a)|=\max_{x\in \R^2}\frac{|\mathcal{R}_a(x)|}{\|\mathcal{R}_a\|_{L^\infty(\R^2)}}=1$. In view of the definition of the operator $\widetilde{\mathcal{L}}_a$, by the maximum principle one can deduce from (\ref{2.24d}) and (\ref{2.24e}) that $|y_a|\le C$ uniformly in $0<a<a^*$.

On the other hand, the usual elliptic regularity theory yields that there exist a subsequence, still denoted by $\{\mathcal{\bar R}_a\}$, of  $\{\mathcal{\bar R}_a\}$ and a function $\mathcal{\bar R}_0\in H^1(\R^2)$ such that $\mathcal{\bar R}_a \to \mathcal{\bar R}_0$ weakly in $H^1(\R^2)$ and strongly in $L^q_{loc}(\R^2)$ for all $q\in [2,\infty)$ as $a\nearrow a^*$. By (\ref{2NN:3}), we obtain from (\ref{2.24dd})--(\ref{2.24e}) that $\mathcal{\bar R}_0$ satisfies
\[
\nabla \mathcal{\bar R}_0(0)=0,\quad \widetilde{\mathcal{L}}\mathcal{\bar R}_0(x)=0  \,\ \mbox{in}\,\ \R^2,
\]
which  yields that $\mathcal{\bar R}_0=\sum ^2_{i=1}c_i\frac{\partial w}{\partial y_i}$ in view of (\ref{2:2:4}).
Since $\nabla \mathcal{\bar R}_0(0)=0$, we obtain that
\[
 \Big(\frac{\partial ^2w(0)}{\partial x_i\partial x_j}\Big)\begin{pmatrix}  c_1\\ c_2 \end{pmatrix} =0.
\]
Because $det\Big(\frac{\partial ^2w(0)}{\partial x_i\partial x_j}\Big)\not =0 $, we deduce from above that $c_1=c_2=0$, and hence $\mathcal{\bar R}_0(x)\equiv 0$ in $\R^2$. This however contradicts to the above fact that  by passing to a subsequence if necessary,  $1=\mathcal{\bar R}_0(y_a)\to \mathcal{\bar R}_0(\bar y_0)$ as $a\nearrow a^*$ for some $\bar y_0\in \R^2$. Therefore, the claim (\ref{2.24c}) holds true.

We next claim that
\begin{equation}\label{2.24}
\big| \mathcal{R}_a(x)\big|,\ \ \big|\nabla \mathcal{R}_a(x)\big|\leq C (\e^4+|\beta_a|) e^{-\frac{1}{11} |x|}\quad \hbox{uniformly in}\, \ \R^2\,\ \hbox{as}\, \ a\nearrow a^*.
\end{equation}
Actually, we deduce from  \eqref{2.24b} and (\ref{2NN:3}) that
\[
\widetilde{\mathcal{L}}_a\frac{\mathcal{R}_a}{ \e^4+|\beta_a|  }=\frac{N_a(x) }{ \e^4+|\beta_a|  }  \leq C_1e^{-\frac {1}{10}|x|}\quad\hbox{uniformly in}\,\ \R^2\,\ \hbox{as}\, \ a\nearrow a^*.
\]
On the other hand, for sufficiently large $R>1$ there exists a constant $C_2=C_2(R)>0$ such that
\begin{equation}\label{2MM:D12}
\frac{\mathcal{R}_a}{ \e^4+|\beta_a|  }\le C_2e^{-\frac {1}{10}|x|} \quad\hbox{at}\,\ |x|=R\,\ \hbox{as}\, \ a\nearrow a^*,
\end{equation}
and
\[
C_2\widetilde{\mathcal{L}}_a e^{-\frac {1}{10}|x|}\ge C_1e^{-\frac {1}{10}|x|}
 \quad\hbox{in}\,\ \R^2/B_R(0)\,\ \hbox{as}\, \ a\nearrow a^*.
\]
We thus have
\begin{equation}\label{2M:D12}
\widetilde{\mathcal{L}}_a  \Big(\frac{\mathcal{R}_a}{ \e^4+|\beta_a|  }- C_2e^{-\frac {1}{10}|x|}\Big)\le 0
 \quad\hbox{in}\,\ \R^2/B_R(0)\,\ \hbox{as}\, \ a\nearrow a^*.
\end{equation}
By the comparison principle, we then obtain from (\ref{2MM:D12}) and (\ref{2M:D12}) that
\[ \frac{\mathcal{R}_a}{ \e^4+|\beta_a|  }\le C_2e^{-\frac {1}{10}|x|} \quad\hbox{in}\,\ \R^2/B_R(0)\,\ \hbox{as}\, \ a\nearrow a^*. \]
Similarly, one can derive that there exists a constant $C_3=C_3(R)>0$ such that
\[
\frac{\mathcal{R}_a}{ \e^4+|\beta_a|  }\ge -C_3e^{-\frac {1}{10}|x|} \quad\hbox{in}\,\ \R^2/B_R(0)\,\ \hbox{as}\, \ a\nearrow a^*,
\]
and hence
\[
|\mathcal{R}_a|\le C(\e^4+|\beta_a|) e^{-\frac {1}{10}|x|} \quad\hbox{in}\,\ \R^2/B_R(0)\,\ \hbox{as}\, \ a\nearrow a^*.
\]
This therefore yields that $\mathcal{R}_a$ satisfies
\begin{equation}\label{2:D12}
|\mathcal{R}_a(x)|\leq C(\e^4+|\beta_a|) e^{-\frac {1}{10}|x|}\quad\hbox{in}\,\ \R^2,
\end{equation}
due to the estimate (\ref{2.24c}). By gradient estimates of (3.15) in \cite{GT}, we further derive from \eqref{2.24b}, (\ref{2NN:3}) and (\ref{2:D12}) that
\begin{equation}\label{2:D8}
|\nabla \mathcal{R}_a|\leq C(\e^4+|\beta_a|) e^{-\frac {1}{11}|x|}\quad\hbox{in}\,\ \R^2.
\end{equation}
Therefore, we obtain from (\ref{2:D12}) and (\ref{2:D8}) that (\ref{2.24}) holds true.


{\em Step 2.}  Since $\psi_{2}(|x|):=-\frac{1}{2}\big(w+x\cdot \nabla w\big)$ is radially symmetric, we deduce from (\ref{2.24R}) that
\[
x^\bot\cdot\nabla R_a=x^\bot\cdot\nabla \big[\mathcal{R}_a +\beta_a\psi_{2}+\e^4\psi_{1}(x) \big]=x^\bot\cdot\nabla [\mathcal{R}_a+\e^4\psi_{1}(x)]
\quad\hbox{in}\ \ \R^2.\]
This then implies from (\ref{2.24}) that
\begin{equation}\label{2:27}
|x^\bot\cdot\nabla R_a|\leq C(\e^4+|\beta_a|)e^{-\frac {1}{12}|x|}\quad \hbox{uniformly in}\, \  \R^2\,\ \hbox{as}\, \ a\nearrow a^*.
\end{equation}
 Similar to (\ref{2:D12}) and (\ref{2:D8}), one can deduce from (\ref{2.21}) and (\ref{2:27}) that
\begin{equation}\label{2:28}
\big| I_a(x)\big|,\ \ \big| \nabla I_a(x)\big| \leq C\e^2(\e^4+|\beta_a|)e^{-\frac {1}{14}|x|}  \quad \hbox{uniformly in}\, \ \R^2\,\ \hbox{as}\, \ a\nearrow a^*,
\end{equation}
where the property (\ref{2:2:2}) is used in view of the fact that $I_a$ satisfies $ \inte wI_adx\equiv 0$.

{\em Step 3.}  Applying (\ref{2:28}), the same argument of (\ref{2NN:3}) yields that
the nonhomogeneous term $N_a(x)$ of (\ref{2.24b}) satisfies
\begin{equation}\label{2NN:5}
\frac{|N_a(x)|}{\e ^4+ |\beta_a|}\le C \delta_a e^{-\frac{1}{10}|x|}\ \ \mbox{uniformly in}\ \ \R^2\,\ \hbox{as}\, \ a\nearrow a^*,
\end{equation}
where  $\delta_a>0$  satisfies $\delta_a=o(1)$ as $a\nearrow a^*$. The same argument of Step 1 then gives from (\ref{2NN:5}) that
\begin{equation}\label{2NN:9}
|\mathcal{R}_a|,\ \ |\nabla \mathcal{R}_a|\leq C(\e^4+|\beta_a|)\delta_a e^{-\frac {1}{12}|x|}\quad\hbox{in}\,\ \R^2,
\end{equation}
where  $\delta_a>0$ also satisfies $\delta_a=o(1)$ as $a\nearrow a^*$. This proves (\ref{2N:a3}) in view of (\ref{2.24R}).

Finally, similar to (\ref{2:27}), one can derive from (\ref{2.24R}) and (\ref{2NN:9}) that
\[|x^\bot\cdot\nabla R_a|\leq C(\e^4+\gamma_a)e^{-\frac {1}{14}|x|}\quad \hbox{uniformly in}\, \  \R^2\,\ \hbox{as}\, \ a\nearrow a^*,
\]
where $\gamma_a>0$ satisfies $\gamma_a=o(|\beta_a|)$ as $a\nearrow a^*$. Applying this refined estimate, the same argument of (\ref{2:28}) then yields that (\ref{2.3:3A}) holds true, and Lemma \ref{lem2.2} is thus proved.\qed

Before ending this subsection, we derive the refined estimate of $I_a$ as $a\nearrow a^*$.

\begin{lem}\label{lem2.3}
Under the assumptions of Proposition \ref{prop2.1}, the imaginary part $I_a$  of (\ref{2:R}) satisfies
\begin{equation}\label{2N:K1}
I_a(x):=\e ^6\,\Omega\,\psi_{I}(x) +o(\e ^6+\e ^2|\beta_a|) \quad \hbox{in}\, \ \R^2 \ \ \mbox{as} \ \ a\nearrow a^*,
\end{equation}
where   $\psi_{I}(x)\in C^2(\R^2)\cap L^\infty(\R^2)$ solves uniquely
\begin{equation}\label{2N:K2}
\inte \psi_{I}wdx=0,\quad  \mathcal{L} \psi_{I}(x)=-  \big(x^\perp \cdot \nabla \psi_1\big)
 \,\ \mbox{in}\,\ \R^2,
\end{equation}
where the operator $\mathcal{L}$ is defined by (\ref{2:2:1}), and $\psi_1(x)\in C^2(\R^2)\cap L^\infty(\R^2)$ is given by (\ref{lem2.1:2}).
\end{lem}

\noindent{\bf Proof.} Following the argument of Lemma \ref{lem2.2}, one can derive from (\ref{2.21}), (\ref{2.24R}) and (\ref{2NN:9}) that $\inte I_awdx\equiv 0$, and
\begin{equation}\label{2N:K3}
\mathcal{L}_aI_a=-\e^2 \Omega\,  \big(x^\perp \cdot \nabla R_a\big)=- \e^6 \Omega\,  \big(x^\perp \cdot \nabla \psi_1\big)+o(\e^6+\e^2|\beta_a|) \ \ \mbox{as} \ \ a\nearrow a^*,
\end{equation}
where $\psi_1(x)\in C^2(\R^2)\cap L^\infty(\R^2)$ is given by (\ref{lem2.1:2}).
Set
\[I _{1,a}(x)=I_a(x)-\e ^6\Omega\,\psi_{I}(x),\]
where $\psi_I(x)\in C^2(\R^2)\cap L^\infty(\R^2)$ is defined in (\ref{2N:K2}).

Similar to the proof of Lemma \ref{lem2.2},    one can deduce from (\ref{2.3:3A}) and (\ref{2N:K3}) that
$I _{1,a}(x)=o(\e ^6+\e ^2|\beta_a|)$ uniformly in $\R^2$ as $a\nearrow a^*$, which then implies directly that (\ref{2N:K1}) holds true.
Finally, since $\mathcal{L}$ is a linear operator, the uniqueness of $\psi_{I}(x)\in C^2(\R^2)\cap L^\infty(\R^2)$ defined by (\ref{2N:K2}) follows from the constriction $\inte \psi_{I}wdx=0$ and the property (\ref{2:2:2}) (cf. \cite[Lemma 4.1]{Wei96}).
\qed

\begin{rem}\label{rem2.1} Under the general
assumptions on $V_\Omega(x)$ of Proposition \ref{prop2.1}, the upper bound (\ref{2.3:3A}) of $I_a$ as $a\nearrow a^*$ is optimal and $I_a$ does not vanish. However, if the trap $V_\Omega(x)$ has a better symmetry, then the leading term $\psi_{I}(x)$ of $I_a$ as $a\nearrow a^*$ may vanish and hence the upper bound (\ref{2.3:3A}) of $I_a$ as $a\nearrow a^*$ is not optimal. Specially, if the trap $V_\Omega(x)$ is radially symmetric, we have proved in \cite[Theorem 1.3]{GLY} that $I_a(x)\equiv 0$  as $a\nearrow a^*$.
\end{rem}

\section{Asymptotic Expansions and Applications}
In this section we first address the proof  of Theorem  \ref{thm2} on asymptotic expansions of minimizers for $e_F(a)$ as $a\nearrow a^*$, based on which we then establish Theorem  \ref{thm1}. Towards this purpose, throughout this section we always consider the harmonic trap $V(x)$ of the form (\ref{5:V}), and let $0<\Omega<\Omega ^*:=2\Lambda <2$ be fixed, so that all estimates of previous section are applicable.

Let $R_a(x)$ and $I_a(x)$ be defined by (\ref{2:v}) and (\ref{2:R}), respectively, where $u_a$ denotes a complex-valued minimizer of $e_F(a)$ as $a\nearrow a^*$. Applying (\ref{5:2}), it then follows from Lemma \ref{lem2.2} that the real part $R_a$ of (\ref{2:R}) satisfies
\begin{equation}\label{5:4a}
R_a(x):=\e ^4\psi_{1}(x)+\beta _a\psi_{2}(|x|) +o(\e ^4+|\beta _a|)\quad \hbox{in}\, \ \R^2 \,\ \mbox{as} \,\ a\nearrow a^*,
\end{equation}
where $\e >0$ and $\beta _a\in \R$ are defined by (\ref{2:beta}).  Here $\psi_{2}(|x|)\in C^2(\R^2)\cap L^\infty(\R^2)$ satisfying
\begin{equation}\label{5:4b2}
\psi_{2}(|x|):=-\frac{1}{2}\big(w+x\cdot \nabla w\big)
\end{equation}
is radially symmetric, and $\psi_{1}(x)\in C^2(\R^2)\cap L^\infty(\R^2)
$
solves uniquely
\begin{equation}\label{5:4b}
 \nabla \psi_{1}(0)=0,\quad \widetilde{\mathcal{L}}\psi_{1}(x)=- \frac{\lam _0^4}{a^*}w^3(x)-\big(x_1^2+\Lambda ^2x_2^2\big) w(x)
 \,\ \mbox{in}\,\ \R^2,
\end{equation}
where $0<\Lambda <1$ is as in (\ref{5:V}), $\lam _0>0$ is as in (\ref{5:3}) in view of (\ref{5:2}), (\ref{5:2K}) and (\ref{2.e}), and while the operator $\widetilde{\mathcal{L}}$ is defined by (\ref{2:2:3}).

We next denote $\psi_{11}(x)\in C^2(\R^2)\cap L^\infty(\R^2)$ to be the unique solution of
\begin{equation}\label{5:4M}
 \nabla \psi_{11}(0)=0,\quad \widetilde{\mathcal{L}}\psi_{11}(x)=- \frac{\lam _0^4}{a^*}w^3(x)- \Lambda ^2|x|^2 w(x)
 \,\ \mbox{in}\,\ \R^2,
\end{equation}
 and let $\psi (x)\in C^2(\R^2)\cap L^\infty(\R^2)$ be the unique solution of
\begin{equation}\label{5:4N}
 \nabla \psi(0)=0,\quad \widetilde{\mathcal{L}}\psi(x)= - (1-\Lambda ^2)x_1^2 w(x)
 \,\ \mbox{in}\,\ \R^2,\ \   x=(x_1,x_2)\in\R^2,
\end{equation}
where $0<\Lambda <1$ is as above. Following (\ref{5:4b}), we decompose $\psi_1$ as
\begin{equation}\label{5:4Q}
\psi_1(x)=\psi_{11}(|x|)+\psi (x), \ \ \mbox{where}\ \ \psi_{11}(|x|)  \ \ \mbox{is radially symmetric in}\, \ \R^2.
\end{equation}
Applying (\ref{5:4Q}), we thus conclude from Lemma \ref{lem2.3} that  the imaginary part $I_a$  of (\ref{2:R}) satisfies
\begin{equation}\label{5:4c}
I_a(x):=\e ^6\,\Omega\,\psi_{I}(x) +o(\e ^6+\e ^2|\beta_a|) \quad \hbox{in}\, \ \R^2 \ \ \mbox{as} \ \ a\nearrow a^*,
\end{equation}
where   $\psi_{I}(x)\in C^2(\R^2)\cap L^\infty(\R^2)$ solves uniquely
\begin{equation}\label{5:4d}
\inte \psi_{I}wdx=0,\quad  \mathcal{L} \psi_{I}(x)=-  \big(x^\perp \cdot \nabla \psi_1\big)=-  \big(x^\perp \cdot \nabla \psi \big)
 \,\ \mbox{in}\,\ \R^2,
\end{equation}
the operator $\mathcal{L}$ is defined by (\ref{2:2:1}), and  $\psi (x)\in C^2(\R^2)\cap L^\infty(\R^2)$ is given uniquely by (\ref{5:4N}).

\begin{lem}\label{lem2.4} Under the assumptions of Theorem \ref{thm2},  the unique maximum point $x_a$ of $|u_a(x)|$ satisfies
\begin{equation}\label{5:5}
\big|\e ^3   x_a \big|= o\big([\e ^4+|\beta _a|]^2\big)\ \ \mbox{as} \ \ a\nearrow a^*,
\end{equation}
where we denote $o\big([\e ^4+|\beta _a|]^2\big)=o\big( \e ^8\big)+o\big( \e ^4 |\beta _a| \big)+o\big( \beta _a^2\big)$.
\end{lem}

\noindent{\bf Proof.} We calculate from (\ref{2.21a}) and (\ref{2.21b}) that
\begin{equation}\label{3a:21}\arraycolsep=1.5pt\begin{array}{lll}
 \displaystyle\inte\frac{\partial w}{\partial x_1}\widetilde{\mathcal{L}}_aR_a&=&\displaystyle\inte\frac{\partial w}{\partial x_1}F_a(x)\\[2mm]
&=&-\e ^4 \Big[\displaystyle\inte\frac{\partial w}{\partial x_1}V_\Omega\big(x+\frac{ x_a}{\e }\big)w+\displaystyle\frac{\lam _0^4}{a^*}\inte\frac{\partial w}{\partial x_1}\widetilde{R}_a^2w\Big]
\\[3mm]
&&+ \eps^4_a\,\Omega \displaystyle\inte\frac{\partial w}{\partial x_1}\Big(x^\bot\cdot\frac{\nabla I_a}{\eps^2_a}\Big)+\frac{a}{a^*}\displaystyle\inte\frac{\partial w}{\partial x_1}I_a^2w\\[3mm]
&:=&A_1+A_2,
\end{array}\end{equation}
where $V_\Omega (x)=V(x)-\frac{\Omega^2}{4}|x|^2$ is as before, the part $A_2$ satisfies
\begin{equation}\label{3a:21b}
\arraycolsep=1.5pt\begin{array}{lll}
A_2 :&=& \eps^4_a\,\Omega  \displaystyle\inte\frac{\partial w}{\partial x_1}\Big(x^\bot\cdot\frac{\nabla I_a}{\eps^2_a}\Big)+\frac{a}{a^*}\displaystyle\inte\frac{\partial w}{\partial x_1}I_a^2w\\[3mm]
&=&\e^4\big[\e^4+o(\e^4+|\beta _a|)\big]\Omega^2  \displaystyle\inte\frac{\partial w}{\partial x_1}\big(x^\bot\cdot \nabla \psi_I\big)\ \ \mbox{as} \ \ a\nearrow a^*
\end{array}\end{equation}
in view of Lemmas \ref{lem2.2} and \ref{lem2.3}, and $\psi_{I}(x)\in C^2(\R^2)\cap L^\infty(\R^2)$ satisfies  (\ref{5:4d}) after simplification. As for the part $A_1$, we observe from (\ref{5:2}) and Proposition \ref{prop2.1} that
$$\inte\frac{\partial w^2}{\partial x_1}V_\Omega(x)dx=0,\quad \inte\frac{\partial w}{\partial x_1}w^2 \psi_{2}(|x|)dx=0,$$
where the radial symmetry of $\psi_{2}(|x|)$ is also used. We then have
\begin{equation}\label{3a:21c}
\arraycolsep=1.5pt\begin{array}{lll}
-\displaystyle\frac{1}{\e^4}A_1:&=&\displaystyle\inte\frac{\partial w }{\partial x_1}V_\Omega\Big(x+\frac{ x_a}{\e }\Big)w+\displaystyle\frac{\lam_0^4}{a^*}\inte\frac{\partial w}{\partial x_1}\widetilde{R}_a^2w\\[3mm]
&=&\displaystyle\inte\frac{\partial w }{\partial x_1}\Big[V_\Omega\Big(x+\frac{ x_a}{\e }\Big)-V_\Omega(x)\Big]w+\displaystyle\frac{\lam_0^4}{a^*}\inte\frac{\partial w }{\partial x_1}\big(\widetilde{R}_a^2-w^2\big)w\\[3mm]
&=&\displaystyle\inte\frac{\partial w }{\partial x_1}\Big(\frac{ x_a}{\e }\cdot\nabla V_\Omega(x)\Big)w+\displaystyle\frac{2\lam_0^4}{ a^*}\e^4\inte\frac{\partial w}{\partial x_1} w^2\psi_1\\[3mm]
&&
+o\Big(\Big|\displaystyle\frac{ x_a}{\e }\Big|+\e^4+|\beta _a|\Big)\ \ \mbox{as} \ \ a\nearrow a^* ,
\end{array}\end{equation}
where the expansion of Lemma \ref{lem2.2}(2) is used in the last identity. We deduce from (\ref{3a:21})--(\ref{3a:21c}) that
\begin{equation}\label{3a:21d}
\arraycolsep=1.5pt\begin{array}{lll}
&& -\e^4\displaystyle\inte\frac{\partial w }{\partial x_1}\Big(\frac{ x_a}{\e }\cdot\nabla V_\Omega(x)\Big)w\\[3mm]
&=&\displaystyle\e^8\Big\{\frac{2\lam_0^4}{ a^*}\inte\frac{\partial w}{\partial x_1} w^2\psi_1-\Omega^2  \displaystyle\inte\frac{\partial w}{\partial x_1}\big(x^\bot\cdot \nabla \psi_I\big)\Big\}\\[3mm]
&&+\displaystyle\inte\frac{\partial w}{\partial x_1}\widetilde{\mathcal{L}}_aR_a+o(\e^8+\e^4|\beta_a|)
\ \ \mbox{as} \ \ a\nearrow a^* ,
\end{array}\end{equation}
where  $\psi_{I}(x)\in C^2(\R^2)\cap L^\infty(\R^2)$ satisfies (\ref{5:4d}).

On the other hand, by the definition of $\widetilde{\mathcal{L}},$ we have $\inte\frac{\partial w}{\partial x_1}\widetilde{\mathcal{L}}R_adx=0$. It then follows from Lemmas \ref{lem2.2} and \ref{lem2.3} that
\[\arraycolsep=1.5pt
\begin{array}{lll}
\displaystyle\inte\frac{\partial w}{\partial x_1}\widetilde{\mathcal{L}}_aR_a&=& \displaystyle\inte\frac{\partial w}{\partial x_1}\big(\widetilde{\mathcal{L}}_a-\widetilde{\mathcal{L}}\big)R_a \\[3mm]
&=& \displaystyle\inte\frac{\partial w}{\partial x_1}\Big\{\e ^4\Big(\displaystyle\frac{\Omega^2}{4}|x|^2+\displaystyle\frac{\lam _0^4}{a^*}w^2\Big)+\e ^4V_\Omega \big(x+\frac{x_a}{\eps_a}\big)-\beta _a\\[3mm]
&&\qquad \qquad\displaystyle -wR_a-\frac{a}{a^*}(2w+R_a)R_a-\frac{a}{a^*}I_a^2\Big\}R_a\\[3mm]
&=& \e ^8 \displaystyle\inte\frac{\partial w}{\partial x_1}\Big(\displaystyle\frac{\Omega^2}{4}|x|^2+\displaystyle\frac{\lam_0 ^4}{a^*}w^2\Big)\psi_1+\e ^8\displaystyle\inte\frac{\partial w}{\partial x_1}V_\Omega(x)\psi_1\\[3mm]
&&+\e ^4\beta_a\displaystyle\inte\frac{\partial w}{\partial x_1}V_\Omega(x)\psi_2-\e ^4\beta_a\displaystyle\inte\frac{\partial w}{\partial x_1} \psi_1-\frac{\e^8}{2}\displaystyle\inte\frac{\partial w^2}{\partial x_1} \psi_1^2\\[3mm]
&&-\e ^4\beta_a\displaystyle\inte\frac{\partial w^2}{\partial x_1} \psi_1\psi_2-\displaystyle\frac{a}{a^*}\inte\frac{\partial w^2}{\partial x_1}\big(\e^8\psi_1^2+2\e ^4\beta_a \psi_1\psi_2\big)\\[3mm]
&&+o([\e^4+|\beta_a|]^2)\ \ \mbox{as} \ \ a\nearrow a^* ,
\end{array}\]
due to the radial symmetry of $\psi_{2}(|x|)$. The above estimate thus gives that
\begin{equation}\label{3a:22}
\arraycolsep=1.5pt
\begin{array}{lll}
\displaystyle\inte\frac{\partial w}{\partial x_1}\widetilde{\mathcal{L}}_aR_a&=& \e ^8\Big\{\displaystyle\inte\frac{\partial w}{\partial x_1}\Big[\displaystyle\frac{\lam _0 ^4}{a^*}w^2+\Big(\displaystyle\frac{\Omega^2}{4}|x|^2+V_\Omega(x)\Big)\Big]\psi_1-\frac{3}{2}\displaystyle\inte\frac{\partial w^2}{\partial x_1} \psi_1^2 \Big\}\\[3mm]
&&+\e ^4\beta_a\Big\{\displaystyle\inte\frac{\partial w}{\partial x_1}V_\Omega(x)\psi_2-\displaystyle\inte\frac{\partial w}{\partial x_1} \psi_1-3\displaystyle\inte\frac{\partial w^2}{\partial x_1} \psi_1\psi_2\Big\}\\[4mm]
&&+o([\e^4+|\beta_a|]^2)\quad \mbox{as} \ \ a\nearrow a^* .
\end{array}
\end{equation}
Applying (\ref{5:2K}), we now derive from (\ref{3a:21d}) and (\ref{3a:22}) that
\begin{equation}\label{5:6}
\arraycolsep=1.5pt\begin{array}{lll}
&& -\displaystyle\frac{\e^4}{2}\inte\frac{\partial w ^2}{\partial x_1}\Big(\frac{ x_a}{\e }\cdot\nabla V_\Omega(x)\Big)\\[3mm]
&=& \e ^8\Big\{\displaystyle\inte\frac{\partial w}{\partial x_1}\Big[\displaystyle\frac{3\lam _0^4}{a^*}w^2+\big(x_1^2+\Lambda ^2x_2^2\big)\Big]\psi_1\\[3mm]
&&\qquad \quad -3\displaystyle\inte\frac{\partial w }{\partial x_1}w \psi_1^2 -\Omega^2  \displaystyle\inte\frac{\partial w}{\partial x_1}\big(x^\bot\cdot \nabla \psi_I\big)\Big\}\\[3mm]
&&-\e ^4\beta_a\Big\{\displaystyle\inte\frac{\partial w}{\partial x_1} \psi_1+6\displaystyle\inte\frac{\partial w }{\partial x_1} w\psi_1\psi_2-\displaystyle\inte\frac{\partial w}{\partial x_1}V_\Omega(x)\psi_2\Big\}\\[4mm]
&&+o([\e^4+|\beta_a|]^2)
\quad \mbox{as} \ \ a\nearrow a^* ,
\end{array}\end{equation}
where  $\psi_{I}(x)\in C^2(\R^2)\cap L^\infty(\R^2)$ is given by (\ref{5:4d}).



Similar to (3.23) of \cite{GLW}, one can obtain that
\begin{equation}\label{3a:27a}
\displaystyle\inte\frac{\partial w}{\partial x_1} \psi_1+6\displaystyle\inte\frac{\partial w }{\partial x_1} w\psi_1\psi_2-\displaystyle\inte\frac{\partial w}{\partial x_1}V_\Omega(x)\psi_2=0.
\end{equation}
Moreover, one can note from (\ref{5:4b}) that $\psi_{1}(x)$ is even in $x\in\R^2$, which thus implies that
\[
\displaystyle\inte\frac{\partial w}{\partial x_1}\Big[\displaystyle\frac{3\lam _0^4}{a^*}w^2+\big(x_1^2+\Lambda ^2x_2^2\big)\Big]\psi_1 -3\displaystyle\inte\frac{\partial w }{\partial x_1}w \psi_1^2=0.
\]
Applying (\ref{5:6}), we then derive from above that
the unique maximum point $x_a$ of $|u_a(x)|$ satisfies
\begin{equation}\label{5:7}
\arraycolsep=1.5pt\begin{array}{lll}
&&\displaystyle\frac{\e^3x_a}{2}\cdot \frac{\partial\nabla H_\Omega (y)}{\partial y_1}  \Big | _{y=0}=
\displaystyle\frac{\e^4}{2}\inte w ^2 (x)\Big(\frac{ x_a}{\e }\Big)\cdot\frac{\partial \nabla V_\Omega(x+y)}{\partial x_1}\Big | _{y=0}dx
\\[3mm] &=& -\displaystyle\frac{\e^4}{2}\inte\frac{\partial w ^2}{\partial x_1} \Big(\frac{ x_a}{\e }\Big)\cdot\nabla V_\Omega(x)dx\\[3mm]
 & =&- \e ^8 \,\Omega^2  \displaystyle\inte\frac{\partial w}{\partial x_1}\big(x^\bot\cdot \nabla \psi_I\big)dx
 +o([\e^4+|\beta_a|]^2)
\quad \mbox{as} \ \ a\nearrow a^* ,
\end{array}\end{equation}
due to the definition of $ H_\Omega (y):=\inte w^2(x)V_\Omega (x+y)dx$.
Similarly, one can obtain that the unique maximum point $x_a$ of $|u_a(x)|$ also satisfies
\begin{equation}\label{5:72}
\arraycolsep=1.5pt\begin{array}{lll}
&&\displaystyle\frac{\e^3x_a}{2}\cdot \frac{\partial\nabla H_\Omega (y)}{\partial y_2}  \Big | _{y=0}=\displaystyle\frac{\e^4}{2}\inte w ^2 (x)\Big(\frac{ x_a}{\e }\Big)\cdot\frac{\partial \nabla V_\Omega(x+y)}{\partial x_2}\Big | _{y=0}dx\\[3mm]
&=& -\displaystyle\frac{\e^4}{2}\inte\frac{\partial w ^2}{\partial x_2} \Big(\frac{ x_a}{\e }\Big)\cdot\nabla V_\Omega(x)dx\\[3mm]
 & =&- \e ^8 \,\Omega^2  \displaystyle\inte\frac{\partial w}{\partial x_2}\big(x^\bot\cdot \nabla \psi_I\big)
dx+o([\e^4+|\beta_a|]^2)
\quad \mbox{as} \ \ a\nearrow a^* .
\end{array}\end{equation}
In the appendix, we shall prove that the following claim is true:
\begin{equation}\label{5:8}
II_1=\inte\frac{\partial w}{\partial x_1}\big(x^\bot\cdot \nabla \psi_I\big)
 =0, \ \ II_2=\inte\frac{\partial w}{\partial x_2}\big(x^\bot\cdot \nabla \psi_I\big)
=0,
\end{equation}
where $\psi_{I}(x)\in C^2(\R^2)\cap L^\infty(\R^2)$ solves uniquely (\ref{5:4d}). Recall from (\ref{5:2}) that $0$ is the unique and non-degenerate critical point of $ H_\Omega (y)$. Applying (\ref{5:8}), we then conclude from (\ref{5:7}) and (\ref{5:72}) that (\ref{5:5}) holds true, and the proof of Lemma \ref{lem2.4} is therefore complete.
\qed

Since the limit estimates of Lemma \ref{lem2.2} are not enough for establishing Theorem \ref{thm2}, we next employ Lemma \ref{lem2.4} to derive the following more terms of $R_a$ in terms of $\e$ and $\beta_a$.

\begin{lem}\label{lem3.1}
Under the assumptions of Theorem \ref{thm2},
the real part $R_a$ of (\ref{2:R}) satisfies
\begin{equation}\label{5:2:1}\arraycolsep=1.5pt
\begin{array}{lll}
R_a(x):&=&\e ^4\psi_{1}(x)+\beta _a\psi_{2}(x) +\e ^8\psi_{3}(x)+\beta _a^2\psi_{4}(x)+\e ^4 \beta _a\psi_{5}(x)\\[2mm]
&&+o\big([\e ^4+|\beta _a|]^2\big)\quad \hbox{in}\, \ \R^2  \ \ \mbox{as} \ \ a\nearrow a^*,
\end{array}\end{equation}
where $\psi_{1}(x)$ and  $\psi_{2}(x) $ are as in (\ref{5:4b2}) and (\ref{5:4b}), respectively. However, $\psi_{i}(x)\in C^2(\R^2)\cap L^\infty(\R^2)$ solves uniquely
\begin{equation}\label{5:2:2}
 \nabla \psi_{i}(0)=0,\quad \widetilde{\mathcal{L}}\psi_{i}(x)=f_i(x)
 \,\ \mbox{in}\,\ \R^2,\ i=3,\,4,\,5,
\end{equation}
and $f_i(x)$ satisfies
\begin{equation}\label{5:2:3}\arraycolsep=1.5pt
f_i(x)=\left\{\begin{array}{lll}
&3w\psi^2_1-\Big[\displaystyle\frac{3\lam _0^4}{a^*}w^2+\big(x_1^2+\Lambda ^2x_2^2\big)\Big]\psi_1+\Omega\, \big(x^\bot\cdot \nabla \psi_{I}\big),\qquad \quad \ &\mbox{if}\ \ i=3;\\[3mm]
 &\psi_2+3w\psi^2_2, \,\ &\mbox{if}\ \ i=4;\\[2mm]
 &  6w\psi_1\psi_2+\psi _1-\Big[\displaystyle\frac{3\lam _0^4}{a^*}w^2+\big(x_1^2+\Lambda ^2x_2^2\big)\Big]\psi_2  , \,\ &\mbox{if}\ \ i=5;
\end{array}\right.\end{equation}
where $\psi_{I}(x)\in C^2(\R^2)\cap L^\infty(\R^2)$ is given by (\ref{5:4d}).
\end{lem}

\noindent{\bf Proof.} Denote \[
M_a:=R_a-\e ^4\psi_1-\beta_a\psi_2,
\]
which then yields from (\ref{2.21a}) that
\begin{equation}\label{2.3:4}
\arraycolsep=1.5pt\begin{array}{lll}
 \widetilde{\mathcal{L}}_aM_a&=& \widetilde{\mathcal{L}}_aR_a- \widetilde{\mathcal{L}}_a(\e ^4\psi_1+\beta_a\psi_2)\\[2mm]
 &=&\big\{\widetilde{\mathcal{L}}_aR_a-\widetilde{\mathcal{L}}(\e ^4\psi_1+\beta_a\psi_2)\big\}-\big(\widetilde{\mathcal{L}}_a-\widetilde{\mathcal{L}}\big)(\e ^4\psi_1+\beta_a\psi_2)\\[2mm]
 &:=&II_1+II_2.
\end{array}\end{equation}
Direct calculations give that the term $II_2$ of (\ref{2.3:4}) satisfies
\begin{equation}\label{2.3:7}
\arraycolsep=1.5pt\begin{array}{lll}
II_2:&=&-\big(\widetilde{\mathcal{L}}_a-\widetilde{\mathcal{L}}\big)(\e ^4\psi_1+\beta_a\psi_2)\\[2mm]
&=&-(\e ^4\psi_1+\beta_a\psi_2)\Big\{\e ^4\Big[\displaystyle\frac{\Omega^2}{4}|x|^2+V_\Omega \big(x+\frac{x_a}{\eps_a}\big)\Big]-\beta _a-wR_a\\[3mm]
&&\qquad\qquad\qquad\quad\ +\displaystyle\frac{\alp_a}{a^*}w^2-\frac{a}{a^*}(2w+R_a)R_a-\frac{a}{a^*}I_a^2\Big\}.
\end{array} \end{equation}
Applying Lemmas \ref{lem2.2} and \ref{lem2.3}, we thus get that
\begin{equation}\label{2.3:7B}
\arraycolsep=1.5pt\begin{array}{lll}
II_2:&=& \e ^8\displaystyle\Big\{3w\psi^2_1-\Big(\frac{\lam  _0^4}{a^*}w^2+\big(x_1^2+\Lambda ^2x_2^2\big)\Big)\psi_1\Big\}\\[2mm]
 &&+\e ^4\beta_a\displaystyle\Big\{ 6w\psi_1\psi_2+\psi _1-\Big(\frac{\lam  _0^4}{a^*}w^2+\big(x_1^2+\Lambda ^2x_2^2\big)\Big)\psi_2\Big\}\\[3mm]
 &&+\beta_a^2\big(\psi_2+3w\psi^2_2\big)+o([\e ^4+|\beta _a|]^2) \ \ \mbox{as} \ \ a\nearrow a^*.
\end{array} \end{equation}

As for the term $II_1$ of (\ref{2.3:4}), we derive from (\ref{2.21a}) that
\begin{equation}\label{2.3:5}
\arraycolsep=1.5pt\begin{array}{lll}
II_1:&=&\widetilde{\mathcal{L}}_aR_a-\widetilde{\mathcal{L}}(\e ^4\psi_1+\beta_a\psi_2)\\[2mm]
&=&\widetilde{\mathcal{L}}_aR_a-\beta _aw+\e ^4\Big\{\displaystyle\frac{\lam  _0^4}{a^*}w^3+\big(x_1^2+\Lambda ^2x_2^2\big)w\Big\}\\[2mm]
&=&-\e ^4\Big[V_\Omega \big(x+\displaystyle\frac{x_a}{\eps_a}\big)-V_\Omega (x)\Big]w-\displaystyle\frac{\alp_a}{a^*}w(2w+R_a)R_a\\[3mm]
&&+\eps^4_a\,\Omega\, \displaystyle\Big(x^\bot\cdot\frac{\nabla I_a}{\eps^2_a}\Big)+\displaystyle\frac{a}{a^*}I_a^2w.
\end{array} \end{equation}
Applying (\ref{5:4c}) yields that
\[
\e^4\,\Omega\, \displaystyle\Big(x^\bot\cdot\frac{\nabla I_a}{\eps^2_a}\Big)= \e ^8\,\Omega\, \big(x^\bot\cdot \nabla \psi_{I}\big)+o(\e^8+\e^4|\beta _a|)\quad \mbox{as} \ \ a\nearrow a^*,
\]
where $\psi_{I}(x)\in C^2(\R^2)\cap L^\infty(\R^2)$ is given by (\ref{5:4d}). By Lemma \ref{lem2.4}, we have
\begin{equation}\label{2.3:6}\arraycolsep=1.5pt\begin{array}{lll}
-\e ^4\displaystyle\Big[V_\Omega \big(x+\frac{x_a}{\eps_a}\big)-V_\Omega (x)\Big]w
&=&-\e ^3\big(x_a \cdot\nabla V_\Omega (x)\big)w(x)\big[1+o(1)\big]\\
&=&o\big([\e ^4+|\beta _a|]^2\big)   \ \ \mbox{as} \ \ a\nearrow a^*.
\end{array} \end{equation}
Applying Lemmas \ref{lem2.2} and \ref{lem2.3}, we now deduce from above that
\begin{equation}\label{2.3:8B}
\arraycolsep=1.5pt\begin{array}{lll}
II_1:&=& \e ^8\displaystyle\Big\{-\frac{2\lam _0^4}{a^*}w^2\psi_1+\Omega\, \big(x^\bot\cdot \nabla \psi_{I}\big)\Big\}\\[3mm]
&&-\e ^4\beta_a\displaystyle\frac{2\lam _0^4}{a^*}w^2\psi_2+o\big([\e ^4+|\beta _a|]^2\big)   \ \ \mbox{as} \ \ a\nearrow a^*,
\end{array} \end{equation}
where $\psi_{I}(x)\in C^2(\R^2)\cap L^\infty(\R^2)$ is given by (\ref{5:4d}).

Applying (\ref{2.3:7B}) and (\ref{2.3:8B}), we now calculate from (\ref{2.3:4}) that
\begin{equation}\label{3a:33}\arraycolsep=1.5pt\begin{array}{lll}
\widetilde{\mathcal{L}}_aM_a&=&\e ^8\displaystyle\Big\{3w\psi^2_1-\Big[\frac{3\lam _0^4}{a^*}w^2+\big(x_1^2+\Lambda ^2x_2^2\big)\Big]\psi_1 +\Omega\, \big(x^\bot\cdot \nabla \psi_{I}\big)\Big\}\\[2mm]
 &&+\e ^4\beta_a\displaystyle\Big\{ 6w\psi_1\psi_2+\psi _1-\Big(\frac{3\lam _0^4}{a^*}w^2+\big(x_1^2+\Lambda ^2x_2^2\big)\Big)\psi_2\Big\}\\[3mm]
 &&+\beta_a^2\big(\psi_2+3w\psi^2_2\big)+o([\e ^4+|\beta _a|]^2) \ \ \mbox{as} \ \ a\nearrow a^*.
\end{array}
\end{equation}
Following (\ref{3a:33}), the argument of Lemma \ref{lem2.2} then yields the estimate (\ref{5:2:1}). Moreover, the property (\ref{2:2:4}) implies the uniqueness of $\psi_i$ for $i=3,4,5$, and the proof of Lemma \ref{lem3.1} is therefore complete.
\qed

\begin{lem}\label{lem3.2}
Under the assumptions of Theorem \ref{thm2}, we have
\begin{equation}\label{2.4:1}
 \displaystyle\inte w\psi_{1} =0,\ \ \inte w\psi_{2} =0,\ \ T_1= \inte \big(2w\psi_4+\psi^2_2\big)=0,
\end{equation} and
 \begin{equation}\label{2.4:3}
T_2=2\inte w\psi_5+2\inte \psi_1\psi_2=-2\lam  _0^4<0,
\end{equation}
where $\psi_1(x), \cdots ,\psi_5(x),\,\psi_I(x)\in C^2(\R^2)\cap L^\infty(\R^2)$ are as in Lemma \ref{lem3.1}.
\end{lem}

\noindent\textbf{Proof.} The assumptions of Theorem \ref{thm2} imply that
\begin{equation}\label{2.4:3K}
\inte  w^2\Big(x\cdot\nabla V(x)\Big)=2 \inte   V(x)w^2dx=2\lam _0^4>0.
\end{equation}
Applying (\ref{2.4:3K}), Lemma \ref{lem3.2} can be proved in a similar way of \cite[Lemma 3.5]{GLW}, and the detailed proof is omitted for simplicity.
\qed

\subsection{Proofs of Theorems \ref{thm2} and \ref{thm1}}

In this subsection we complete the proofs of Theorems \ref{thm2} and \ref{thm1}. Under the assumption (\ref{5:V}), we remark from (\ref{5:2K}) that $\e =\ep>0$, where $\ep>0$ is as in the statement of Theorem \ref{thm2}.

 \vskip 0.05truein

\noindent\textbf{Proof  of Theorem \ref{thm2}.}   Applying Proposition \ref{prop2.1}, we derive from (\ref{2:R}) that $R_a$ satisfies
\begin{equation}\label{2a:k}
 \inte w^2=\inte  |v_a|^2=\inte \big[\big(w+R_a\big)^2+I^2_a\big],\  i.e., \ \ 0=2\inte wR_a+\inte R_a^2+\inte I_a^2.
\end{equation}
Following Lemmas \ref{lem3.1} and \ref{lem3.2}, we then derive from (\ref{2a:k}) that
\begin{equation}\label{3a:37}\arraycolsep=1.5pt\begin{array}{lll}
0&=&2\displaystyle\inte wR_a+\displaystyle\inte R_a^2+\displaystyle\inte I_a^2\\[3mm]
&=&2\displaystyle\inte w(\e ^4\psi_{1}  +\beta _a\psi_{2}+ \e ^8\psi_{3} +\beta ^2_a\psi_{4} +\e ^4\beta _a\psi_{5})\\[4mm]
&&+\displaystyle\inte (\e ^4\psi_{1} +\beta _a\psi_{2} + \e ^8\psi_{3} +\beta ^2_a\psi_{4} +\e ^4\beta _a\psi_{5})^2+o([\e ^4+|\beta _a|]^2)\\[4mm]
&=& 2\e ^4\Big(\displaystyle \inte  w\psi_{1}\Big)+2\beta _a\Big(\displaystyle  \inte  w\psi_{2}\Big)+\beta _a^2\Big(2\displaystyle\inte  w\psi_{4}+\displaystyle\inte \psi_{2}^2\Big)\\[4mm]
&&+2\e ^4\beta _a\Big(\displaystyle\inte w\psi_{5}+\displaystyle\inte \psi_{1}\psi_{2}\Big)+\e ^8\Big(2\displaystyle\inte w\psi_{3}+\displaystyle\inte \psi_{1}^2\Big)+o([\e ^4+|\beta _a|]^2)\\[4mm]
&=&-2\lam _0^4\e ^4\beta _a+\e ^8\Big(2\displaystyle\inte w\psi_{3}+\displaystyle\inte \psi_{1}^2\Big)+o([\e ^4+|\beta _a|]^2) \ \ \mbox{as} \ \ a\nearrow a^*,
\end{array}
\end{equation}
where (\ref{5:4c}) and (\ref{5:4d}) are also used. One can derive from (\ref{3a:37}) that
\[
2\displaystyle\inte w\psi_{3}+\displaystyle\inte \psi_{1}^2\not=0,
\]
and
\begin{equation}\label{5:2:beta1}
-\displaystyle2\lam _0^4 \beta _a+\e ^4 \Big(2\displaystyle\inte w\psi_{3}+\displaystyle\inte \psi_{1}^2\Big)=0,
\end{equation}
where $\psi_{1}(x)$ and  $\psi_{3}(x) $ are as in (\ref{5:4b}) and (\ref{5:2:2}), respectively.
In the appendix, we shall prove the following claim that
\begin{equation}\label{5:2:beta}
I:=2\displaystyle\inte w\psi_{3}+\displaystyle\inte \psi_{1}^2=\displaystyle\inte (3w^2-1)\psi^2_1- \displaystyle 4\inte \big(x_1^2+\Lambda ^2x_2^2\big)w \psi_1.
\end{equation}
Together with (\ref{5:2:beta}), we conclude from (\ref{5:2:beta1}) that the constant $\beta_a$ satisfies
\begin{equation}\label{5:2:beta2}
\beta _a=C^*\e ^4,\ \ \mbox{where}\ \ C^*=\frac{1 }{2\lam _0^4}\Big[\displaystyle\inte (3w^2-1)\psi^2_1- \displaystyle 4\inte \big(x_1^2+\Lambda ^2x_2^2\big)w \psi_1\Big]\not =0.
\end{equation}
Applying (\ref{5:2:1}), (\ref{5:4c}) and (\ref{5:2:beta2}), we thus conclude from (\ref{2:v}) and (\ref{2:R}) that the refined limit profile (\ref{thm2:3}) holds true. Moreover, the estimate (\ref{thm2:4}) follows directly from (\ref{5:2:beta2}) and Lemma \ref{lem2.4}. This completes the proof of Theorem \ref{thm2}.\qed

As a byproduct of Theorem \ref{thm2}, we finally establish Theorem \ref{thm1} on the nonexistence of vortices in a very large region.

\vskip 0.05truein

\noindent\textbf{Proof  of Theorem \ref{thm1}.}  By Theorem \ref{thm2}, for any fixed and sufficiently large $R>1$, there exists a constant $C_R:=C(R)>0$ such that
\begin{equation}\label{3F:2M}
|v_{a}(x)|\ge w(x)-C_R\eps_a ^4>0\ \ \mbox{in}\ \ \big\{x\in\R^2:\, |x|\le R\big\}  \ \ \mbox{as} \ \ a\nearrow  a^*,
\end{equation}
where $v_{a}(x)$ is as in (\ref{2:v}).

Setting $\tilde{w}_{a}=v_{a}(x)-w(x)$, we claim that
\begin{equation}\label{OK:9}
    |\tilde{w}_{a}|
    \leq C_1 \eps_{a}^{4}|x|^{\frac{5}{2}}e^{-\sqrt{1-C_2\eps_{a}^{4}}|x|}\ \ \hbox{uniformly in}\ \ \R^2 \ \ \mbox{as} \ \ a\nearrow  a^*,
\end{equation}
where the constants $C_1>0$ and $C_2>0$ are independent of $0<a<a^*$.
To prove the claim (\ref{OK:9}), we note from (\ref{2:va}) that $\tilde{w}_{a}$ satisfies
\begin{equation*}
    (-\triangle+\hat{V}_{a})(\tilde{w}_{a}+w)+i\eps_{a}^{2}\Omega (x^{\bot}\cdot\nabla\tilde{w}_{a})=0
    \ \ \hbox{in}\ \ \R^2,
\end{equation*}
i.e.,
\begin{equation}\label{1}
    (-\triangle+\hat{V}_{a})\tilde{w}_{a}+i\eps_{a}^{2}\Omega (x^{\bot}\cdot\nabla\tilde{w}_{a})+(-\triangle+\hat{V}_{a})w=0\ \ \hbox{in}\ \ \R^2,
\end{equation}
where
\begin{equation}\label{2}
    \hat{V}_{a}(x)=\frac{\eps_{a}^{4}\Omega^{2}}{4}|x|^2
    +\eps_{a}^{2}V_{\Omega}(\eps_{a}x+x_{a})-\mu_{a}\eps_{a}^{2}-\frac{a}{a^*}|v_{a}|^{2}\ \ \hbox{in}\ \ \R^2.
\end{equation}
One can derive from (\ref{1}) that
\begin{equation}\label{3}
\begin{split}
    &-\frac{1}{2}\triangle |\tilde{w}_{a}|^{2}+\Big[\frac{\eps_{a}^{4}\Omega^{2}}{4}|x|^2
    +\eps_{a}^{2}V_{\Omega}(\eps_{a}x+x_{a})-\mu_{a}\eps_{a}^{2}-\frac{a}{a^*}|v_{a}|^{2}\Big]|\tilde{w}_{a}|^{2}
    +|\nabla\tilde{w}_{a}|^{2}\\
    &-\eps_{a}^{2}\Omega x^{\bot}(i\tilde{w}_{a},\nabla\tilde{w}_{a})+(-\triangle+\hat{V}_{a})(w,\tilde{w}_{a})=0\ \ \hbox{in}\ \ \R^2.
\end{split}
\end{equation}
By the diamagnetic inequality (\ref{1:2:1A}), we have
\begin{equation}\label{4}
    |\nabla\tilde{w}_{a}|^{2}+\frac{\eps_{a}^{4}\Omega^{2}}{4}|x|^2|\tilde{w}_{a}|^{2}
    -\eps_{a}^{2}\Omega x^{\bot}(i\tilde{w}_{a},\nabla\tilde{w}_{a})
    \geq \big|\nabla|\tilde{w}_{a}|\big|^{2}\ \ \hbox{in}\ \ \R^2.
\end{equation}
Since
\[
 \frac{1}{2}\triangle |\tilde{w}_{a}|^{2}=|\tilde{w}_{a}|\Delta |\tilde{w}_{a}|+\big|\nabla |\tilde{w}_{a}|\big|^2,
\]
we deduce from (\ref{3}) and (\ref{4}) that
\begin{equation}\label{5}
    -\triangle|\tilde{w}_{a}|-\mu_{a}\eps_{a}^{2}|\tilde{w}_{a}|
    \leq |(-\triangle+\hat{V}_{a})w|+\frac{a}{a^*}|v_{a}|^{2}|\tilde{w}_{a}|\ \ \hbox{in}\ \ \R^2.
\end{equation}

The argument of \cite[Proposition 3.3]{GLY} gives that as $a\nearrow  a^*,$
\[|v_{a}|^{2}\le Ce^{-\frac{4}{3}|x|}\ \ \hbox{in}\ \ \R^2,\]
which implies from Theorem \ref{thm2} that as $a\nearrow  a^*,$
\begin{equation}\label{6}
    \frac{a}{a^*}|v_{a}|^{2}|\tilde{w}_{a}|\leq C\eps_{a}^{4}e^{-\frac{4}{3}|x|}\ \ \hbox{in}\ \ \R^2.
\end{equation}
Note from (\ref{2:beta}) and (\ref{5:2:beta2}) that as $a\nearrow  a^*,$
\[
 -\mu_{a}\eps_{a}^{2}=1-C^*\eps_{a}^{4}.
\]
We then calculate from (\ref{2:exp}) and (\ref{2}) that  as $a\nearrow  a^*,$
\begin{equation}\label{7}
\begin{split}
   |(-\triangle+\hat{V}_{a})w|&=\Big|\Big(\frac{\eps_{a}^{4}\Omega^{2}}{4}|x|^2
    +\eps_{a}^{2}V_{\Omega}(\eps_{a}x+x_{a})-1-\mu_{a}\eps_{a}^{2}+w^{2}-\frac{a}{a^*}|v_{a}|^{2}\Big)w\Big|\\
    &\leq C \eps_{a}^{4}|x|^{\frac{3}{2}}e^{-|x|}\ \, \ \hbox{in}\ \ \R^2\backslash B_R(0),
\end{split}
\end{equation}
where the sufficiently large constant $R>1$ is as in (\ref{3F:2M}).
We thus deduce from (\ref{5})--(\ref{7}) that  as $a\nearrow  a^*,$
\begin{equation}\label{8}
    -\triangle|\tilde{w}_{a}|-\mu_{a}\eps_{a}^{2}|\tilde{w}_{a}|
    \leq C_0\eps_{a}^{4}|x|^{\frac{3}{2}}e^{-|x|}\ \, \ \hbox{in}\ \ \R^2\backslash B_R(0),
\end{equation}
where the constant $C_0>0$ is independent of $0<a<a^*$.
Since Theorem \ref{thm2} gives that $|\tilde{w}_{a}|=O(\eps_{a}^{4})$ as $a\nearrow  a^*,$ we have  as $a\nearrow  a^*,$
\begin{equation}\label{8M}
|\tilde{w}_{a}| \le C  \eps_{a}^{4}|x|^{\frac{5}{2}}e^{-\sqrt{1-|C^*|\eps_{a}^{4}}\,|x|}\ \ \hbox{at}\ \ |x|=R>1,
\end{equation}
where $R>1$ is as in (\ref{3F:2M}), and $C>0$ is large enough and independent of $0<a<a^*$.
By the comparison principle as in (\ref{2:D12}), we thus derive from (\ref{8}) and (\ref{8M}) that the claim (\ref{OK:9}) holds true, in view of the fact that $|\tilde{w}_{a}|=O(\eps_{a}^{4})$ as $a\nearrow  a^*$.

Applying (\ref{2:exp}) and (\ref{OK:9}), we now have as  $a\nearrow  a^*,$
\begin{equation}\label{10}
\begin{split}
    &|v_{a}|\geq |w|-|v_{a}-w|\\
     \geq &C_{3}|x|^{-\frac{1}{2}}e^{-|x|}-C_{1}\eps_{a}^{4}|x|^{\frac{5}{2}}e^{-\sqrt{1-C_2\eps_{a}^{4}}\,|x|}\\
    \ge &|x|^{-\frac{1}{2}}e^{-|x|}\big(C_3-C_1\eps_{a}^{4}|x|^{3}e^{C\eps_{a}^{4}|x|}\big)>0, \ \ \mbox{if}\ \ R\le|x|\le \Big(\frac{C_3}{2 C_1}\Big)^\frac{1}{3}\e^{-\frac{4}{3}}.
\end{split}
\end{equation}
We thus conclude from (\ref{3F:2M}) and (\ref{10}) that as  $a\nearrow  a^*,$
\[
|v_a(x)|=\Big|\eps_{a}\sqrt{a^*}u_{a}(\eps_{a}x+x_{a})e^{-i(\frac{\eps_{a}\Omega}{2}x\cdot x_{a}^{\bot}-\theta_{a})}\Big|>0\ \ \mbox{in}\ \ \Big\{x\in\R^2:\, |x|\le \Big(\frac{C_3}{2 C_1}\Big)^\frac{1}{3}\e^{-\frac{4}{3}}\Big\}.
\]
Thus, there exists a small constant  $C_*>0$, independent of $0<a<a^*$, such that  as  $a\nearrow  a^*,$
\begin{equation}\label{11}
    |u_{a}(y)|>0 \ \ \hbox{if}\ \  |y|\leq C_{4}\epsilon_{a}^{-\frac{1}{3}}\leq \frac{C_*}{(a^*-a)^{\frac{1}{12}}},
\end{equation}
i.e., $u_a$ does not admit any vortex in the region $R(a):=\big\{x\in\R^2:\, |x|\le C_*(a^*-a)^{-\frac{1}{12}} \big\}$ as $a\nearrow a^*$. This completes the proof of Theorem \ref{thm1}.
\qed

\appendix
\section{Appendix}

In this appendix, we follow those notations of Section 3 to address the proof of the  claims (\ref{5:8}) and (\ref{5:2:beta}).
\vskip 0.05truein

\noindent{\bf Proof of   (\ref{5:8}).} Consider the polar coordinate $(r, \theta)$ in $ \R^2$, where $\theta \in [0, 2\pi ]$.
Rewrite $w(x)=w(r)$, $\psi (x)=\psi (r, \theta)$ and $\psi_{I}(x)=\psi_{I}(r, \theta)$, where $\psi $ and $\psi_{I} $ are given by (\ref{5:4N}) and (\ref{5:4d}), respectively.  We then have
\begin{equation}\label{5:9a}
\frac{\partial w}{\partial x_1}=w'(r) \cos \theta, \ \ \frac{\partial w}{\partial x_2}=w'(r)\sin \theta.
\end{equation}
Since $\nabla \psi =\frac{x}{r} \psi  _r+\frac{x^\bot}{r^2} \psi _\theta$, we have
\begin{equation}\label{5:9b}
x^\bot\cdot \nabla \psi=\frac{\partial \psi (r,\theta)}{\partial\theta} , \ \ x^\bot\cdot \nabla \psi_I=\frac{\partial \psi_I (r,\theta)}{\partial\theta}.
\end{equation}
By the symmetry of the linear inhomogeneous equation (\ref{5:4N}), we deduce that
\begin{equation}\label{5:10}
\psi (r, \theta)=\psi (r, 2\pi -\theta), \ \ \psi (r, \theta)=\psi (r,  \theta-\pi),\ \ \theta \in [\pi, 2\pi ],
\end{equation}
where $\psi (x)\in C^2(\R^2)\cap L^\infty(\R^2)$ is given by (\ref{5:4N}) as before. It  then yields from (\ref{5:10}) that
\begin{equation}\label{5:11}
\frac{\partial \psi (r,\theta)}{\partial\theta}=-\frac{\partial \psi (r, 2\pi -\theta)}{\partial (2\pi -\theta)}, \ \ \frac{\partial \psi (r,\theta)}{\partial\theta}=\frac{\partial\psi (r,  \theta-\pi)}{\partial(\theta-\pi)},\ \ \theta \in [\pi, 2\pi ].
\end{equation}
Note from (\ref{5:4d}) and (\ref{5:9b}) that
$\psi_{I}(x)\in C^2(\R^2)\cap L^\infty(\R^2)$ is the unique solution of
\begin{equation}\label{5:12}
\inte \psi_{I}wdx=0,\quad  \mathcal{L} \psi_{I}(x)=-  \big(x^\perp \cdot \nabla \psi \big)=-\frac{\partial \psi (r,\theta)}{\partial\theta}
 \,\ \mbox{in}\,\ \R^2.
\end{equation}
We then derive from (\ref{5:11}) and (\ref{5:12}) that the unique solution $\psi_{I}(x) $ satisfies
\begin{equation}\label{5:12a}
\psi _I (r, \theta)=-\psi_I (r, 2\pi -\theta), \ \ \psi _I(r, \theta)=\psi _I(r,  \theta-\pi),\ \ \theta \in [\pi, 2\pi ].
\end{equation}

Applying (\ref{5:9b}),  we now deduce   that
\begin{equation}\label{5:12E}\arraycolsep=1.5pt\begin{array}{lll}
II_1&=&\displaystyle\inte\frac{\partial w}{\partial x_1}\big(x^\bot\cdot \nabla \psi_I\big)\\[3mm]
&=&\displaystyle\int ^\infty _0\int ^{2\pi}_0rw'(r)\cos \theta \frac{\partial \psi _I(r,\theta)}{\partial\theta}d\theta dr\\[3mm]
&=&\displaystyle\int ^\infty _0\int ^{2\pi}_0rw'(r)\sin \theta   \psi _I(r,\theta) d\theta dr\\[3mm]
&=&\displaystyle \int ^\infty _0\int ^{\pi}_0rw'(r)\sin \theta   \psi _I(r,\theta) d\theta dr+\displaystyle\int ^\infty _0\int ^{2\pi}_{\pi}rw'(r)\sin \theta   \psi _I(r,\theta) d\theta dr\\[3mm]
&:=&A_1+A_2.
\end{array}\end{equation}
As for the term $A_2$, we derive from (\ref{5:12a}) that
\begin{equation}\label{5:12F}\arraycolsep=1.5pt\begin{array}{lll}
A_2:&=&\displaystyle\int ^\infty _0\int ^{2\pi}_{\pi}rw'(r)\sin \theta   \psi _I(r,\theta) d\theta dr\\[3mm]
&=&-\displaystyle\int ^\infty _0\int ^{2\pi}_{\pi}rw'(r)\sin (\theta -\pi)   \psi _I(r,\theta -\pi) d \theta  dr\\[3mm]
&=&\displaystyle-\int ^\infty _0\int ^{\pi}_{0}rw'(r)\sin \delta  \, \psi _I(r,\delta) d \delta  dr=-A_1,
\end{array}\end{equation}
where we denote $\delta :=\theta -\pi$. We hence obtain from (\ref{5:12E}) and (\ref{5:12F}) that
\begin{equation}\label{5:15A}
II_1=\displaystyle\inte\frac{\partial w}{\partial x_1}\big(x^\bot\cdot \nabla \psi_I\big)=A_1+A_2=0.
\end{equation}

Similar to (\ref{5:12E}), one can conclude from (\ref{5:9b}) and (\ref{5:12a}) that
\begin{equation}\label{5:15B}\arraycolsep=1.5pt\begin{array}{lll}
II_2&=&\displaystyle\inte\frac{\partial w}{\partial x_2}\big(x^\bot\cdot \nabla \psi_I\big)\\[3mm]
&=&-\displaystyle \int ^\infty _0\int ^{\pi}_0rw'(r)\cos \theta   \psi _I(r,\theta) d\theta dr-\displaystyle\int ^\infty _0\int ^{2\pi}_{\pi}rw'(r)\cos \theta   \psi _I(r,\theta) d\theta dr\\[3mm]
&:=&-(B_1+B_2),
\end{array}\end{equation}
where the term $B_2$ satisfies
\begin{equation}\label{5:12C}\arraycolsep=1.5pt\begin{array}{lll}
B_2:&=&\displaystyle\int ^\infty _0\int ^{2\pi}_{\pi}rw'(r)\cos \theta   \psi _I(r,\theta) d\theta dr\\[3mm]
&=&-\displaystyle\int ^\infty _0\int ^{2\pi}_{\pi}rw'(r)\cos (2\pi-\theta  )   \psi _I(r,2\pi-\theta ) d \theta  dr\\[3mm]
&=&\displaystyle\int ^\infty _0\int _{\pi}^{0}rw'(r)\cos \delta  \, \psi _I(r,\delta) d \delta  dr=-B_1,
\end{array}\end{equation}
where we denote $\delta :=2\pi -\theta$. The above estimates yield that
\[
II_2=\displaystyle\inte\frac{\partial w}{\partial x_2}\big(x^\bot\cdot \nabla \psi_I\big)=-(B_1+B_2)=0,
\]
together with (\ref{5:15A}), which thus implies that (\ref{5:8}) holds true, and we are done.
\qed

\vskip 0.05truein
In the rest of this appendix, we establish the claim  (\ref{5:2:beta}) as follows.

\vskip 0.05truein

\noindent{\bf Proof of (\ref{5:2:beta}).}
Since $\psi_2$ is radially symmetric, we first note that for $\psi '_2=\frac{d\psi _2}{dr}$,
\[\arraycolsep=1.5pt
\begin{array}{lll}
\displaystyle\Omega\inte  \psi_2\big(x^\bot\cdot \nabla \psi_{I}\big)&=&\displaystyle\Omega\inte \psi_2\big[-x_2(\psi_{I})_{x_1}+x_1(\psi_{I})_{x_2}\big]\\[3mm]
&=&\displaystyle\Omega\inte \big[\psi_Ix_2(\psi_2)_{x_1}-\psi_Ix_1(\psi_{2})_{x_2}\big]\\[3mm]
&=&\displaystyle\Omega\inte \psi_I\frac{x_1x_2}{r}\big(\psi '_2-\psi '_2\big)=0.
\end{array}\]
Following (\ref{5:4b}) and (\ref{5:2:2}), we then derive from above that for $V(x)=x_1^2+\Lambda ^2x_2^2$,
\begin{equation}\label{5:beta:1}
\arraycolsep=1.5pt
\begin{array}{lll}
I&=&2\displaystyle\inte  \psi_{3}\widetilde{\mathcal{L}}\psi_{2}+\displaystyle\inte \psi_{1}^2\\[3mm]
&=&2\displaystyle\inte  \psi_{2}\widetilde{\mathcal{L}}\psi_{3}-2\displaystyle\Omega\inte  \psi_2\big(x^\bot\cdot \nabla \psi_{I}\big)+\displaystyle\inte \psi_{1}^2\\[3mm]
 &=&-\displaystyle\inte  \big(w+x\cdot \nabla w\big)\Big\{3w\psi^2_1-\Big[\displaystyle\frac{3\lam _0^4}{a^*}w^2+V(x)\Big]\psi_1\Big\}+\displaystyle\inte \psi_{1}^2\\[3mm]
 &:=&A+B,
\end{array}\end{equation}
where the part $A$ satisfies
\begin{equation}\label{5:beta:2a}
\arraycolsep=1.5pt
\begin{array}{lll}
A&=&-\displaystyle\inte   w \Big\{3w\psi^2_1-\Big[\displaystyle\frac{3\lam _0^4}{a^*}w^2+V(x)\Big]\psi_1\Big\}+\displaystyle\inte \psi_{1}^2\\[3mm]
&=&\displaystyle-\inte | \nabla \psi_1|^2+\frac{2\lam _0^4}{a^*}\displaystyle\inte w^3\psi_1.
\end{array}\end{equation}
Since
\[
\inte  (x\cdot \nabla \psi_1)\Delta \psi_1 =-2\inte  | \nabla \psi_1|^2-\inte  (x\cdot \nabla \psi_1)\Delta \psi_1, \]
$\ i.e., \
\ -\inte  | \nabla \psi_1|^2 =\inte  (x\cdot \nabla \psi_1)\Delta \psi_1,$
it follows from (\ref{5:beta:2a}) that the part $A$ can be rewritten as
\begin{equation}\label{5:beta:2b}
A=\displaystyle\inte  (x\cdot \nabla \psi_1)\Delta \psi_1+\frac{2\lam _0^4}{a^*}\displaystyle\inte w^3\psi_1.
\end{equation}

We rewrite  the part $B$ of (\ref{5:beta:1}) as
\begin{equation}\label{5:beta:2}
B=-\displaystyle\inte  (x\cdot \nabla w) \Big\{3w\psi^2_1-\Big[\displaystyle\frac{3\lam _0^4}{a^*}w^2+V(x)\Big]\psi_1\Big\}
=B_1+B_2+B_3,
\end{equation}
where the term $B_2$ satisfies
\[
\arraycolsep=1.5pt
\begin{array}{lll}
B_2&=&\displaystyle\frac{3\lam _0^4}{a^*}\inte  (x\cdot \nabla w)  w^2\psi_1= - \displaystyle\frac{2\lam _0^4}{a^*}\inte w^3\psi_1-\displaystyle\frac{\lam _0^4}{a^*}\inte w^3 (x\cdot \nabla \psi_1).
\end{array}\]
Together with (\ref{5:beta:2b}), we obtain from above that
\begin{equation}\label{5:beta:4b}
\arraycolsep=1.5pt
\begin{array}{lll}
A+B_2&=&\displaystyle\inte  (x\cdot \nabla \psi_1)\Delta \psi_1-\displaystyle\frac{\lam _0^4}{a^*}\inte w^3 (x\cdot \nabla \psi_1).
\end{array}\end{equation}
The term $B_1$ of (\ref{5:beta:2}) satisfies
\begin{equation}\label{5:beta:3}
\arraycolsep=1.5pt
\begin{array}{lll}
B_1&=&-\displaystyle\inte  (x\cdot \nabla w) 3w\psi^2_1=3\displaystyle\inte w^2\psi^2_1+3\displaystyle\inte w^2\psi_1(x\cdot \nabla \psi_1)\\[3mm]
&=&  \displaystyle\inte (3w^2-1)\psi^2_1+\displaystyle\inte (3w^2-1)\psi_1(x\cdot \nabla \psi_1),
\end{array}\end{equation}
due to the fact that
\[
-\inte \psi_1(x\cdot \nabla \psi_1)=2\inte \psi_1^2+\inte \psi_1(x\cdot \nabla \psi_1).
\]
But the term $B_3$ of (\ref{5:beta:2}) satisfies
\begin{equation}\label{5:beta:5}
\arraycolsep=1.5pt
\begin{array}{lll}
B_3&=&\displaystyle \inte  \psi_1V(x)(x\cdot \nabla w) \\[3mm]
&=& - \displaystyle 2\inte V(x)w \psi_1- \displaystyle \inte w\psi_1[x\cdot \nabla V(x)]-\displaystyle \inte wV(x)(x\cdot\nabla \psi_1)\\[3mm]
&=&- \displaystyle 4\inte V(x)w \psi_1-\displaystyle \inte wV(x)(x\cdot\nabla \psi_1),
\end{array}\end{equation}
since $V(x)$ satisfies $x\cdot \nabla V(x)=2V(x)$.
Applying (\ref{5:beta:2})--(\ref{5:beta:5}), we now obtain from (\ref{5:4b}) and (\ref{5:beta:1}) that
\begin{equation}\label{5:beta:6}
\arraycolsep=1.5pt
\begin{array}{lll}
I&=&(A+B_2)+B_1+B_3  \\[3mm]
&=&\displaystyle\inte (3w^2-1)\psi^2_1- \displaystyle 4\inte V(x)w \psi_1   \\[3mm]
&&+ \displaystyle\inte (x\cdot\nabla \psi_1)\Big\{(\Delta -1+3w^2)\psi_1-\displaystyle\frac{\lam _0^4}{a^*}  w^3-V(x)w\Big\}
\\[3mm]
&=&\displaystyle\inte (3w^2-1)\psi^2_1- \displaystyle 4\inte V(x)w \psi_1,
\end{array}\end{equation}
and the claim (\ref{5:2:beta}) is therefore proved in view of (\ref{5:V}).\qed

\vskip 0.16truein
\noindent {\bf Acknowledgements:} The author is very grateful to the referees for many valuable suggestions which lead to the great improvements of the present paper. The author also thanks Dr. Yong Luo very much for his fruitful discussions on the present paper.


\bigskip




\begin{thebibliography}{GNN}
	
\bibitem{Abo} J. R. Abo-Shaeer, C. Raman, J. M. Vogels and W. Ketterle, {\em Observation of vortex lattices in Bose-Einstein condensate}, Science {\bf 292} (2001), 476.






\bibitem{A} A. Aftalion,  Vortices in Bose-Einstein condensates, Progress in Nonlinear Differential Equations and their Applications, 67. Birkh$\ddot{a}$user Boston, Inc., Boston, MA,  2006.

\bibitem{AA} A. Aftalion, S. Alama and L. Bronsard, {\em Giant vortex and breakdown of strong pinning in
a rotating Bose-Einstein condenstate}, Arch. Ration. Mech. Anal. {\bf 178} (2005), 247--286.


\bibitem{AJ} A. Aftalion, R. L. Jerrard and J. Royo-Letelier, {\em Non-existence of vortices in the small density region of a condensate}, J. Funct. Anal. {\bf 260} (2011), 2387--2406.

\bibitem{AN} A. Aftalion, B. Noris, and C. Sourdis, {\em Thomas-Fermi approximation for coexisting two component Bose-Einstein condensates and nonexistence of vortices for small rotation}, Comm. Math. Phys.  {\bf 336} (2015), no. 2, 509--579.



\bibitem{Anderson} M. H. Anderson, J. R. Ensher, M. R. Matthews, C. E. Wieman and E. A. Cornell, {\em Observation of Bose-Einstein condensation in a dilute atomic vapor}, Science {\bf 269} (1995), 198--201.



\bibitem{ANS} J. Arbunich, I. Nenciu and C. Sparber, {\em Stability and instability properties of rotating Bose-Einstein condensates}, Lett. Math. Phys. {\bf 109} (2019), 1415--1432.



\bibitem{AS} G. Arioli and A. Szulkin, {\em A semilinear Schr$\ddot{o}$dinger equation in the presence of a magnetic field}, Arch. Ration. Mech. Anal. {\bf 170} (2003), 277--295.



\bibitem{BC} W. Bao and Y. Cai, {\em Ground states of two-component Bose-Einstein condensates with an internal atomic Josephson junction}, East Asia J. Appl. Math. {\bf 1} (2011), 49--81.

\bibitem{BJ} T. Bartsch, L. Jeanjean and N. Soave, {\em Normalized solutions for a system of coupled cubic Schr$\ddot{o}$dinger equations on $\R^3$}, J. Math. Pures Appl. {\bf 106} (2016), no. 4, 583--614.

\bibitem{BL} T. Bartsch, Y. Y. Liu and Z. L. Liu, {\em Normalized solutions for a class of nonlinear Choquard equations}, SN Partial Differ. Equ. Appl. (2020), 1: 34.

\bibitem{BW} T. Bartsch, X. Zhong and W. Zou, {\em Normalized solutions for a coupled Schr$\ddot{o}$dinger system}, Math. Ann. {\bf 380}  (2021), 1713--1740.


\bibitem{BH}  N. Basharat, H. Hajaiej, Y. Hu and S. J. Zheng, {\em Threshold for blowup and stability for nonlinear Schrodinger equation with rotation}, submitted (2020), 37 pages, https://arxiv.org/abs/2002.04722

\bibitem{BB} F. Bethuel, H. Brezis, and F. Helein, Ginzburg-Landau Vortices, Progress in Nonlinear
Differential Equations and their Applications {\bf 13}, Birkhauser Boston, Inc., Boston, MA, 1994.


\bibitem{B} I. Bloch, J. Dalibard and W. Zwerger, {\em Many-body physics with ultracold gases}, Reviews of Modern Phys. {\bf 80} (2008), 885--964.






\bibitem{Hulet1} C. C. Bradley, C. A. Sackett, J. J. Tollett and R. G. Hulet, {\it Evidence of Bose-Einstein condensation in an atomic gas with attractive interactions}, Phys. Rev. Lett. {\bf 75} (1995), 1687. {\it Erratum} Phys. Rev. Lett. {\bf 79} (1997), 1170.

\bibitem{Cao} D. M. Cao and  Z. W. Tang, {\em Existence and uniqueness of multi-bump bound states of nonlinear Schr\"{o}dinger  equations with electromagnetic fields}, J. Differential Equations  {\bf 222} (2006), no. 2, 381--424.



\bibitem{CC} L. D. Carr and C. W. Clark, {\em Vortices in attractive Bose-Einstein condensates in two dimensions}, Phys.
Rev. Lett. {\bf 97} (2006), 010403.


\bibitem{CD} Y. Castin and R. Dum, {\em Bose-Einstein condensates with vortices in rotating traps}, European
Phys. J. D {\bf 7} (1999), 399--412.


\bibitem{C} T. Cazenave, Semilinear Schr$\ddot{o}$dinger equations, Courant Lecture Notes in Mathematics  Vol. 10, Courant Institute of Mathematical Science/AMS, New York, 2003.

\bibitem{CL} T. Cazenave and P. L. Lions, \textit{Orbital stability of standing waves for some nonlinear Schr\"{o}dinger equations}, Comm. Math. Phys.  {\bf 85}  (1982), no. 4, 549--561.



\bibitem{CP} M. Correggi, F. Pinsker, N. Rougerie, and J. Yngvason, {\em Rotating superfluids in anharmonic
traps: From vortex lattices to giant vortices}, Phys. Rev. A  {\bf 84} (2011), p. 053614.



\bibitem{CR} M. Correggi and N. Rougerie, {\em Boundary behavior of the Ginzburg-Landau order parameter in the surface superconductivity regime}, Arch. Ration. Mech. Anal. {\bf 219} (2016), 553--606.




\bibitem{D} F. Dalfovo, S. Giorgini, L. P. Pitaevskii and S. Stringari, {\em Theory of Bose-Einstein condensation in trapped gases}, Reviews of Modern Phys. {\bf 71} (1999), 463--512.



\bibitem{EL} M. J. Esteban and P. L. Lions, {\em Stationary solutions of nonlinear Schr$\ddot{o}$dinger equations with an external magnetic field}, Partial differential equations and the calculus of variations, Vol. I,  401--449, Progr. Nonlinear Differential Equations Appl. 1, Birkhuser Boston, Boston, MA, 1989.

\bibitem{F} A. L. Fetter, {\em Rotating trapped Bose-Einstein condensates}, Reviews of Modern Phys. {\bf 81} (2009), 647--691.

\bibitem{Frank} R. L. Frank, {\em Ground states of semi-linear PDEs}, Lecture notes from ``Summer school
on Current Topics in Mathematical Physics", CIRM Marseille, 2013.

\bibitem{GNN} B. Gidas, W. Ni and L. Nirenberg, {\em Symmetry of positive solutions of nonlinear elliptic equations in $\R^n$}, Mathematical analysis and applications  Part A, Adv. in Math. Suppl. Stud. Vol. {\bf 7}, Academic Press, New York  (1981), 369--402.

\bibitem{GT} D. Gilbarg and N. S. Trudinger, Elliptic Partial Differential Equations of Second Order, Springer, 1997.

\bibitem{Gou}   T. X. Gou, Z. T. Zhang, {\em Normalized solutions to the Chern-Simons-Schr\"odinger system},  J. Funct. Anal. {\bf 280} (2021), 108894, 65 pp.




\bibitem{GLW} Y. J. Guo, C. Lin and J. C. Wei, {\em Local uniqueness and refined spike profiles of ground states for two-dimensional attractive Bose-Einstein condensates}, SIAM J. Math. Anal. {\bf 49} (2017), 3671--3715.


\bibitem{GLP} Y. J. Guo, Y. Luo and S. J. Peng, {\em Local uniqueness of ground states for rotating Bose-Einstein condensates with attractive interactions},  Calc. Var. Partial Differential Equations  {\bf 60}  (2021), Paper No. 237, 27pp.

\bibitem{LP} Y. J. Guo, Y. Luo and S. J. Peng, {\em Existence and asymptotic behavior of ground states of rotating Bose-Einstein condensates}, submitted, (2021), 31 pages, arxiv.org/abs/2106.14369

\bibitem{GLY} Y. J. Guo, Y. Luo and W. Yang, {\em The nonexistence of vortices for rotating Bose-Einstein condensates with attractive interactions}, Arch. Rational Mech. Anal. {\bf 238} (2020), 1231--1281.






\bibitem{GS}  Y. J. Guo and R. Seiringer, {\em On the mass concentration for Bose-Einstein condensates with attractive interactions}, Lett. Math. Phys. {\bf 104} (2014), 141--156.



\bibitem{GWZZ} Y. J. Guo, Z. Q. Wang, X. Y. Zeng and H. S. Zhou,  {\em  Properties of ground states of attractive Gross-Pitaevskii equations with multi-well potentials}, Nonlinearity {\bf 31}  (2018),  957--979.







\bibitem{HM}  C. Huepe, S. Metens, G. Dewel, P. Borckmans and  M.E. Brachet, {\em Decay rates in attractive Bose-Einstein condensates}, Phys. Rev. Lett. {\bf 82}  (1999), 1616--1619.



\bibitem{IM-1} R. Ignat and V. Millot, {\em The critical velocity for vortex existence in a two-dimensional rotating Bose-Einstein condensate}, J. Funct. Anal. {\bf 233} (2006), 260--306.



\bibitem{IM-2}   R. Ignat and V. Millot, {\em Energy expansion and vortex location for a two-dimensional rotating Bose-Einstein condensate}, Rev. Math. Phys. {\bf 18} (2006), 119--162.


\bibitem{J}  L. Jeanjean, {\em Existence of solutions with prescribed norm for semilinear elliptic equations}, Nonlinear Anal. {\bf 28} (1997), no. 10, 1633--1659.






%












%
%
%
\bibitem{K} M. K. Kwong, {\em Uniqueness of positive solutions of $\Delta u-u+u^p=0$  in $\R^N$}, Arch. Rational Mech. Anal. {\bf 105} (1989), 243--266.

\bibitem{LNR} M. Lewin, P. T. Nam and N. Rougerie, {\em The mean-field approximation and the nonlinear
Schr\"{o}dinger functional for trapped Bose gases}, Trans. Amer. Math. Soc.   {\bf 368} (2016), 6131--6157.

\bibitem{Lewin} M. Lewin, P. T. Nam and N. Rougerie, {\em Blow-up profile of rotating 2D focusing Bose gases}, Macroscopic Limits of Quantum Systems, Springer Verlag, 2018.





\bibitem{Lieb} E. H. Lieb and M. Loss, Analysis, Graduate Studies in Mathematics Vol. 14. Amer. Math. Soc., Providence, RI, second edition, 2001.



\bibitem{Lieb06} E. H. Lieb and R. Seiringer, {\em Derivation of the Gross-Pitaevskii equation for rotating Bose gases}, Comm. Math. Phys. {\bf 264} (2006), 505--537.





\bibitem {LSS} E. H. Lieb, R. Seiringer, J. P. Solovej and  J. Yngvason, {\em The mathematics of the Bose gas and its condensation}, Oberwolfach Seminars, {\bf 34} Birkh$\ddot{a}$user Verlag, Basel, 2005.



\bibitem {LSY} E. H. Lieb, R. Seiringer and J. Yngvason, {\em Bosons in a trap: A rigorous derivation of the Gross-Pitaevskii energy functional}, Phys. Rev. A {\bf 61} (2000), 043602.





\bibitem {LC} E. Lundh, A. Collin and K.-A. Suominen, {\em Rotational states of Bose gases with attractive interactions
	in anharmonic traps}, Phys. Rev. Lett. {\bf 92} (2004), 070401.



\bibitem {LPY} P. Luo, S. Peng, J. Wei and S. Yan, {\em  Excited states on Bose-Einstein condensates with attractive interactions}, Calc. Var. Partial Differential Equations {\bf 60}  (2021), Paper No. 155, 30pp.




\bibitem{MC1} K. Madison, F. Chevy, J. Dalibard and W. Wohlleben, {\em Vortex formation in a stirred
Bose-Einstein condensate}, Phys. Rev. Lett. {\bf 84} (2000), 806.

\bibitem{MC2} K. Madison, F. Chevy, J. Dalibard and W. Wohlleben, {\em Vortices in a stirred Bose-Einstein condensate}, J. Mod. Opt. {\bf 47} (2000), 2715--2723.




\bibitem{NR} P. T. Nam and N. Rougerie, {\em Improved stability for 2D attractive Bose gases}, J. Math. Phys. {\bf 61} (2020), 021901.




\bibitem{NT} W. M. Ni and I. Takagi, {\em On the shape of least-energy solutions to a semilinear Neumann problem}, Comm. Pure Appl. Math. {\bf 44} (1991), 819--851.


\bibitem{PP} B. Pellacci, A. Pistoia, G. Vaira and G. Verzini, {\em Normalized concentrating solutions to nonlinear elliptic
problems},  J. Differential Equations {\bf 275} (2021), 882--919.


\bibitem{R} N. Rougerie, {\em Non linear Schr\"{o}dinger limit of bosonic ground states, again}, Confluentes Math. {\bf 12} (2020), no. 1, 69--91.

\bibitem{Roug} N. Rougerie, {\em Scaling limits of bosonic ground states, from many-body to nonlinear Schr\"{o}dinger}, EMS Surveys in Math. Sciences  {\bf 7}  (2020), 253--408.


\bibitem{SSbook} E. Sandier and S. Serfaty,  Vortices in the Magnetic Ginzburg-Landau Model, Progress in Nonlinear Differential Equations and their Applications {\bf 70}, Basel: Birkh\'auser, 2007.


\bibitem{S} I. M. Sigal, {\em Magnetic vortices, Abrikosov lattices and automorphic functions},  Math. and comput. modeling, 19--58, Pure Appl. Math. (Hoboken), Wiley, Hoboken, NJ, 2015.




\bibitem{Wei96} J. C. Wei, {\em On the construction of single-peaked solutions to a singularly perturbed semilinear Dirichlet problem}, J. Diff. Eqns. {\bf 129}  (1996),  no. 2, 315--333.




\bibitem{W} M. I. Weinstein, {\em Nonlinear Schr$\ddot{o}$dinger equations and sharp interpolations estimates}, Comm. Math. Phys. {\bf 87} (1983), 567--576.





\bibitem {Z} J. Zhang, {\em Stability of attractive Bose-Einstein condensates}, J. Stat. Phys. {\bf 101} (2000), 731--746.




	
	
\end{thebibliography}
\end{document}